\documentclass[10pt,twocolumn,twoside]{IEEEtran}

\ifCLASSINFOpdf
\else
\fi
\usepackage[latin1]{inputenc}
\usepackage{amssymb}
\usepackage{amscd}
\usepackage{amsmath}
\usepackage{graphicx}
\usepackage{subfigure}
\usepackage{cite}
\usepackage{colortbl}
\usepackage{color}
\usepackage{xcolor}
\usepackage{dsfont}
\usepackage{psfrag}

\newcommand{\nq}{\mathbb{Q}}
\newcommand{\nr}{\mathbb{R}}
\newcommand{\real}{{\cal R}e}
\newcommand{\imag}{{\cal I}m}
\newcommand{\Iz}{\mathbf{I}}
\newcommand{\Az}{\mathbf{A}}
\newcommand{\Bz}{\mathbf{B}}

\newcommand{\bz}{\mathbf{b}}
\newcommand{\Pz}{\mathbf{P}}
\newcommand{\Gz}{\mathbf{G}}
\newcommand{\Cz}{\mathbf{C}}

\newcommand{\xz}{\mathbf{x}}
\newcommand{\qz}{\mathbf{q}}
\newcommand{\wz}{\mathbf{w}}
\newcommand{\Tz}{\boldsymbol{\varUpsilon}}
\newcommand{\tz}{\mathbf{t}}
\newcommand{\pz}{\mathbf{p}}
\newcommand{\uz}{\mathbf{u}}

\newcommand{\vz}{\mathbf{v}}
\newcommand{\Sz}{\mathbf{S}}
\newcommand{\Hz}{\mathbf{H}}
\newcommand{\Zz}{\mathbf{Z}}
\newcommand{\sz}{\mathbf{s}}
\newcommand{\Fz}{\mathbf{F}}
\newcommand{\dz}{\mathbf{d}}
\newcommand{\ez}{\mathbf{e}}
\newcommand{\Wz}{\mathbf{W}}

\newcommand{\etaz}{\boldsymbol{\eta}}
\newcommand{\Lbd}{\boldsymbol{\Lambda}}

\newcommand{\Am}{\boldsymbol{\mathcal{A}}}
\newcommand{\Bm}{\boldsymbol{\mathcal{B}}}
\newcommand{\Cm}{\boldsymbol{\mathcal{C}}}
\newcommand{\Dm}{\boldsymbol{\mathcal{D}}}
\newcommand{\Em}{\boldsymbol{\mathcal{E}}}
\newcommand{\Fm}{\boldsymbol{\mathcal{F}}}
\newcommand{\Gm}{\boldsymbol{\mathcal{G}}}
\newcommand{\Hm}{\boldsymbol{\mathcal{H}}}

\newcommand{\IM}{\boldsymbol{\mathcal{I}}}
\newcommand{\KM}{\boldsymbol{\mathcal{K}}}
\newcommand{\JM}{\boldsymbol{\mathcal{J}}}
\newcommand{\RM}{\boldsymbol{\mathcal{R}}}

\def\0z{{\boldsymbol{0}}}

\title{Low-Complexity Quaternion Adaptive Filters}

\author{Fernando G. Almeida Neto,~\IEEEmembership{Student Member,~IEEE,}
        Vítor H. Nascimento,~\IEEEmembership{Senior Member,~IEEE}
        \thanks{Fernando G. Almeida Neto and Vítor H. Nascimento are with the Department
of Electronic Systems, Escola Politécnica of the University of São Paulo, São Paulo,
SP, Brazil. e-mail: \{fganeto, vitor\}@lps.usp.br.}\thanks{This work is partially funded by CNPq under grants 140994/2011-4 and 245909/2012-5}}

\begin{document}

\maketitle

\begin{abstract}

A general representation of the quaternion gradients presented in the literature is proposed, and an universal update equation
for QLMS-like algorithms is obtained. The general update law is used to study the convergence of widely linear (WL)
algorithms. It is proved that techniques obtained with a gradient similar to the i-gradient are the fastest-converging in 
two situations: 
1) When the correlation matrix contains elements only in 2 axis (1 and $i$, for instance), and
2) When the algorithms use a real data vector, obtained staking up the real and imaginary parts of the original quaternion input vector. 
The general update law is also used to study the convergence of WL-QLMS-based algorithms, and an accurate second-order model is developed for
quaternion algorithms using real-data input. Based on the proposed analysis, we obtain
the fastest-converging WL-QLMS algorithm with real-regressor vector, which is also less costly than the reduced-complexity WL-QLMS (RC-WL-QLMS) algorithm 
proposed in our previous work. It is shown that the new method corresponds to the four-channel LMS algorithm written in the quaternion domain, 
and that they have the same computational complexity. 
Simulations illustrate the performance of the new technique and the accuracy of the analysis.


\end{abstract}

\begin{IEEEkeywords}
    Quaternion processing, quaternion adaptive filtering, widely linear adaptive filtering 
\end{IEEEkeywords}

\section{Introduction}
\label{sec:int} 

Quaternion numbers \cite{kantor1989hypercomplex} were invented by Hamilton in the 19th century, as a generalization of complex numbers to a 
higher-dimensional domain. They consist of one real part and three imaginary elements, usually identified by $i$, $j$ and $k$, where 
$i^2 = j^2 = k^2= -1$. Quaternions appear in many fields, and their applications have been spreading recently, since they can be used to  
concisely describe multi-variable data.

Quaternion algebra is traditionally employed to represent rotations in coordinate systems and to image processing.
In the first application, quaternion algebra provides mathematical robustness to represent rotations. It avoids the gimbal lock 
in Euler angle representations \cite{hanson2006}, which is exploited by attitude control systems \cite{attitude1991, attitude1996}. 
In color image processing, for instance \cite{segmentation2001, color1999}, many techniques employ a quaternion-based model to describe color images, 
allowing a concise representation of the color attributes in a single entity. 
In recent applications, quaternions have been applied to study DNA structures \cite{chang2012dna}, neural networks \cite{mandic_nonlinear2009}, 
beamforming \cite{beamforming} and adaptive filtering \cite{fernandoSSP2011, cctook2008, mandic2010, mandic_recursive2010, took2010quaternion, took2009quaternion}, 
among many others. The last field has experimented a large development lately, and a variety of algorithms have been proposed 
for multi-variable estimation.


There are different forms to define quaternion differentiation (see \cite{Took2011}, \cite{mandic_iQLMS} and \cite{iqlms_review}), 
which gave rise to different
quaternion adaptive algorithms. The first to be proposed was QLMS \cite{took2009quaternion}, whose update equation contains an extra
term (when compared to the complex LMS update law \cite{sayed2008adaptive}) to take into account the non-commutative nature of
quaternion multiplication.
Later, after a new definition of the differentiation to account for quaternion involutions \cite{Took2011}, 
iQLMS \cite{mandic_iQLMS} was proposed and a lower-cost and faster-converging technique emerged. 
However, both QLMS and iQLMS are designed for $\nq$-circular data \cite{santamaria2010}, 
for which the correlation matrix is sufficient to assess full second-order statistics \cite{santamaria2010}.
When the inputs are non $\nq$-circular, these algorithms -- which are generally called strictly linear (SL) --  are not able to 
account for full second-order statistics. In this case, widely linear (WL) algorithms \cite{mandic_recursive2010} can be applied to improve the performance.

Widely linear adaptive filters were initially proposed for complex signal processing \cite{tulay2009}.
In the quaternion case, the definition of WL processes has led to the 
augmented QLMS \cite{took2009quaternion} (which uses the original SL regressor vector and its conjugate as the WL input data) and  to WL-QLMS \cite{took2010quaternion},
where the SL data vector and three quaternion involutions
are the inputs. Later, WL-iQLMS \cite{mandic_asilomar2011} was also proposed as an improved WL-QLMS algorithm.
For all these methods, the WL vector is four times the length of the original SL input vector, and thus the computational cost is 
significantly higher.

In order to reduce the computational complexity of WL-QLMS, we proposed in \cite{fernandoSSP2011} the reduced-complexity (RC) WL-QLMS algorithm. 
The technique was designed to use a real-regressor vector, obtained from the concatenation of the real and imaginary parts of the original quaternion SL data 
vector. With this approach, the algorithm avoids many quaternion-quaternion computations, which are replaced by real-quaternion operations, 
less costly to compute.

In this paper, we generalize our previous work from \cite{fernandoSSP2011} and \cite{fernandoISWCS2012}. We develop an universal description for 
the quaternion gradients proposed in the literature, and we use it to derive the update law of any QLMS-like algorithm.  
The update equation is applied to study the convergence of WL algorithms obtained with different gradients. 
We prove that a class of gradients which includes the i-gradient of \cite{Took2011} leads to the fastest-converging WL-QLMS algorithm, 
under some conditions on the correlation of the input data. We further show that the fastest-converging real-regressor vector WL-QLMS
which is based on the approach of \cite{fernandoISWCS2012} corresponds to the four-channel LMS algorithm (4-Ch-LMS) written in the quaternion 
domain.

The specific contributions of this paper are as follows.
\begin{enumerate}
	\item We propose a general approach to describe the different quaternion gradients proposed in the literature, and we show that different
	      derivatives can generate the same quaternion gradient.
	\item We obtain a general update law, which describes all the QLMS algorithms proposed in the literature. We use it to study the convergence 
	      of WL algorithms and we prove that different derivatives can be used to obtain the same algorithm. We show that gradients 
	      similar to the i-gradient of \cite{Took2011} and to the gradient proposed in \cite{iqlms_review} lead to the fastest-converging WL-QLMS 
	      algorithms in two situations: i) When at most two of the quaternion elements in the input vector are correlated; and
	      ii) When the regressor vector is real and obtained with the concatenation of the real and the 
	      imaginary parts of the original SL quaternion input.
	\item We develop the fastest-converging WL-QLMS algorithm with real-regressor vector, 
	      and we prove that this algorithm corresponds to the 4-Ch-LMS written in the quaternion domain. It is also shown 
	      that the new technique is a reduced-complexity version of WL-iQLMS. 
	\item We show that the proposed algorithm and the 4-Ch-LMS have the same computational complexity, while WL-iQLMS is four times 
	      more costly to compute. 
	\item We extend the second-order analysis of \cite{fernandoISWCS2012} to any WL-QLMS algorithm with a real-regressor vector obtained with 
	      the concatenation of the real and the imaginary parts of the original SL quaternion data vector. The analysis is suitable for 
	      correlated and uncorrelated inputs. 
	      Concise equations 
	      to compute the EMSE (excess mean-square-error) and the MSD (mean square deviation) 
	      \cite{sayed2008adaptive} are also derived. 
	\item We present simulations comparing the performance of the proposed algorithm and other algorithms from the literature. 
	      The second-order model is compared to the algorithms to show the accuracy of our approach.
\end{enumerate}
We note that preliminary results were presented in conference papers (\cite{fernandoSSP2011,fernandoISWCS2012}). 
A reduced-complexity widely-linear complex LMS algorithm was proposed previously in \cite{fernandoISWCS2010}.

The paper is organized as follows. We present a brief review on quaternion algebra and Kronecker products \cite{hornTMA91} in Section 
\ref{sec:preliminaries}, which are applied in our analysis. In Section \ref{sec:estimation}, we review basic concepts of quaternion estimation and 
$\nq$-properness, while Section \ref{sec:quaternion_gradients} introduces our general approach to write quaternion gradients.
We propose our new reduced-complexity algorithm in Section \ref{sec:rc_wl_qlms}, and we develop the analysis of WL quaternion algorithms using 
real regressor data in Sections \ref{sec:mean_analysis} and \ref{sec:mean_square_analysis}. Simulations are presented in Section \ref{sec:simulations}, 
and in Section \ref{sec:conclusions} we conclude the paper.
\\
\\
\textbf{Notation:} We use lower case to describe scalar quantities (e.g.: $a$) and bold lower case to describe column vectors (e.g.: $\mathbf{b}$). 
Bold capital letters represent matrices (e.g.: $\mathbf{A}$). The conjugation of a quaternion is denoted by $(\cdot)^{*}$,  while
$(\cdot)^T$ stands for transposition. $(\cdot)^H$ indicates the transposition and conjugation of a matrix or vector. 
We define the operator $\text{diag}(\cdot)$ for two situations: when the argument is a matrix $\mathbf{A}$, $\text{diag}(\mathbf{A})$ denotes a column vector with the diagonal elements of $\mathbf{A}$.
When the argument is a vector $\mathbf{b}$, $\text{diag}(\mathbf{b})$ denotes a diagonal matrix where the non-zero elements are given by $\mathbf{b}$.
The operations $\imag\{\cdot\}$ and $\real\{\cdot\}$ 
take only the imaginary and the real parts of a complex  number or a quaternion, respectively, and we use the subscripts $\text{R}$, $\text{I}$, $\text{J}$ and $\text{K}$ 
to identify the real and the imaginary parts of a quaternion. $E\{\cdot\}$ is the expectation operator and $\Iz_{N}$ 
is the $N \hspace{-0.11 cm}\times \hspace{-0.11 cm}N$ identity matrix. We use $\text{col}(\cdot)$ to define a column vector.

\section{Preliminaries}
\label{sec:preliminaries}


In this section, we briefly summarize some properties of quaternion algebra and Kronecker products.
These concepts simplify the analysis and the equations derived in this paper.

\subsection{Review on quaternion algebra}
\label{sec:quaternion_algebra}

 A quaternion $q$ is defined as
 \begin{equation*}
   q = q_\text{R} + i q_\text{I} + j q_\text{J} + k q_\text{K},
   \label{eq:quaternion}
 \end{equation*}
where $q_\text{R}$, $q_\text{I}$, $q_\text{J}$ and $q_\text{K}$ are real numbers and $i$, $j$ and $k$ are the imaginary parts, which satisfy $i^2=j^2=k^2=-1$.
The main difference between a quaternion and a complex number is that the multiplication in $\nq$ is not commutative, 
since \cite{kantor1989hypercomplex}
\begin{equation*}
    \begin{array}{ccc}
	  ij = -ji = k,  &  jk = -kj = i, & ki = -ik =j,
    \end{array}
\end{equation*}
such that for two quaternions $q_1$ and $q_2$, in general,
$q_1q_2 \hspace{-0.1 cm} \neq \hspace{-0.1 cm} q_2q_1$.

Similar to complex algebra, the conjugate and the absolute value of a quaternion $q$ are given by
\begin{equation*}
    q^{*} = q_\text{R} - iq_\text{I} - jq_\text{J} - kq_\text{K}
\end{equation*}
and $|q|=\sqrt{qq^{*}}$,
respectively. We can also define the following involutions of $q$, which correspond to\footnote{Note that the conjugate of a quaternion 
$q$ is also an involution of $q$}
\begin{equation*}
    \begin{array}{c}
      q^i \triangleq -iqi = q_\text{R} + i q_\text{I} - j q_\text{J} - k q_\text{K}\\
      q^j \triangleq -jqj = q_\text{R} - i q_\text{I} + j q_\text{J} - k q_\text{K}\\
      q^k \triangleq -kqk = q_\text{R} - i q_\text{I} - j q_\text{J} + k q_\text{K}. 
    \end{array}
\end{equation*}
The involutions appear in the definition of WL algorithms (e.g. \cite{fernandoSSP2011}, 
\cite{took2010quaternion}, \cite{mandic_recursive2010}). They are used to improve the algorithms' performance, when compared to their SL counterparts.

These are the main definitions of quaternion algebra used in this paper. See reference \cite{kantor1989hypercomplex} for a more 
detailed explanation.

\subsection{Properties of Kronecker products}
\label{sec:kronecker}

The Kronecker product \cite{hornTMA91} is an efficient manner to compactly represent some large matrices which have a block-structure.
For the purpose of the analyses performed in this paper, the most relevant properties of Kronecker products -- which are represented by the operator $\otimes$
-- are
\begin{enumerate}
  \item $\mathbf{A} \otimes \left( \mathbf{B} + \mathbf{C} \right) =  \mathbf{A} \otimes \mathbf{B} + \mathbf{A} \otimes  \mathbf{C}$.
  \item $\lambda \left( \mathbf{A} \otimes \mathbf{B}\right) = \left( \lambda \mathbf{A} \right) \otimes \mathbf{B}
				    = \mathbf{A} \otimes \left( \lambda \mathbf{B}\right) $.
  \item $\left( \mathbf{A} \otimes \mathbf{B}\right) \left( \mathbf{C} \otimes \mathbf{D}\right) = \mathbf{A}\mathbf{C} \otimes \mathbf{B}\mathbf{D}$,
	where the number of rows of $\mathbf{C}$ ($\mathbf{B}$) and the number of columns of $\mathbf{A}$ ($\mathbf{D}$) are equal. \label{property_3}
  \item $\left( \mathbf{A} \otimes \mathbf{B}\right)^{-1} = \mathbf{A}^{-1} \otimes \mathbf{B}^{-1}$. \label{property_4}
  \item $\left( \mathbf{A} \otimes \mathbf{B}\right)^{T} = \mathbf{A}^{T} \otimes \mathbf{B}^{T}$.
  \item $\text{Tr}(\mathbf{A} \otimes \mathbf{B})=\text{Tr}(\mathbf{A})\text{Tr}(\mathbf{B})$.\label{property_6}
  \item The eigenvalues of $\left( \mathbf{A} \otimes \mathbf{B}\right)$, where $\mathbf{A}$ is $N \times N$ and $\mathbf{B}$ is $M \times M$, are given by $\lambda_i \eta_j$, for $i=1, \hspace{0.1 cm} 2, \hdots, N$ and 
	$j=1,\hspace{0.1 cm} 2, \hdots, M$, where $\lambda_i$  and $\eta_j$ are the eigenvalues of $\mathbf{A}$ and $\mathbf{B}$, respectively. \label{property_7}
\end{enumerate}

These properties appear implicitly or explicitly in the analyses that follow. They make the equations easier to manipulate and simplify
the interpretation of the resulting expressions.

%
%

\section{Strictly-linear and widely-linear quaternion estimation}
\label{sec:estimation}

In the context of complex WL estimation, the concept of complex circularity is required to 
define when WL algorithms outperform SL ones \cite{picimbono1995}. In this section, we summarize how circularity is extended to 
the quaternion domain, using the comparison between SL and WL quaternion estimation.

Define the $N \times 1$ quaternion data vector
\begin{equation*}
    \qz(n) = \qz_\text{R}(n)+ i\qz_\text{I}(n) + j\qz_\text{J}(n) + k\qz_\text{K}(n).
    \label{eq:q_reg}
\end{equation*}
Given a desired quaternion sequence $d(n)$, the problem solved by a strictly linear estimator is the computation of $\wz_{\text{SL}}(n)$ in  
\begin{equation}
    \hat{d}_{\text{SL}}(n) = \wz_{\text{SL}}^H(n) \qz(n)
    \label{eq:d_SL}
\end{equation}
which minimizes the mean-square error (MSE) $E\{|e_{\text{SL}}(n)|^2\}$, where
\begin{equation}
  e_{\text{SL}}(n) = d(n) - \hat{d}_{\text{SL}}(n).
  \label{eq:e_SL}
\end{equation}
From the orthogonality principle \cite{sayed2008adaptive}, it must be true that
\begin{equation}
    e_{\text{SL}}(n) \perp \qz(n),
    \label{eq:erro_SL2}
\end{equation}
where $\perp$ stands that $e_{\text{SL}}(n)$ and $\qz(n)$ are orthogonal.
Eq. \eqref{eq:erro_SL2} can also be expressed as the expectation
\begin{equation}
    E\{\qz(n) e_{\text{SL}}^*(n)\}=\0z, 
    \label{eq:exp_SL}
\end{equation}
and substituting eq. \eqref{eq:e_SL} in \eqref{eq:exp_SL}, we get
\begin{equation}
    E\{\qz(n) d^*(n)\} = E\{\qz(n) \hat{d}_{\text{SL}}^*(n)\}.
    \label{eq:interm}
\end{equation}
Using eq. \eqref{eq:d_SL} in \eqref{eq:interm}, we can write the system of equations
\begin{equation}
    \Cz_\qz \wz_{\text{SL}}(n) = \pz_\qz,
    \label{eq:system_SL}
\end{equation}
where $\Cz_\qz = E\{\qz(n)\qz^H(n)\}$ is the correlation matrix and $\pz_\qz=E\{\qz(n) d^*(n)\}$ is the cross-correlation vector \cite{picimbono1995}. 

For a quaternion WL estimator, the data vector is modified to account for the involutions, and it is given by\footnote{Note that 
other sets of four different involutions of $\qz(n)$ can be chosen, since they provide similar estimation performance.}
\begin{equation*}
    \qz_{\text{WL}}(n) = \text{col}(\hspace{-0.1 cm}\begin{array}{cccc} \qz(n), & \qz^*(n), & \qz^{i*}(n) & \qz^{j*}(n)\end{array}\hspace{-0.1 cm}),
\end{equation*}
which is four times the length of the original vector $\qz(n)$.
For this approach, one must find the vector $\wz_{\text{WL}}(n)$
which minimizes the MSE condition $E\{|e_{\text{WL}}(n)|^2\}$, where
\begin{align}
    &e_{\text{WL}}(n) = d(n)-\hat{d}_{\text{WL}}(n),
    \label{eq:e_WL}\\
    &\hat{d}_{\text{WL}}(n) = \wz_{\text{WL}}^H(n)\qz_{\text{WL}}(n).
    \label{eq:d_WL}
\end{align}
Again, the orthogonality condition implies that $e_{\text{WL}}(n)$ must be orthogonal to all involutions, i.e., 
\begin{equation}
    E\{\qz_{\text{WL}}(n) e_{\text{WL}}^*(n)\}= \0z.
    \label{eq:exp_WL}
\end{equation}
Using eqs. \eqref{eq:e_WL}, \eqref{eq:d_WL} and \eqref{eq:exp_WL}, we obtain
\begin{equation}
    \Cz_{\text{WL}} \wz_{\text{WL}}(n) = \pz_{\text{WL}},
    \label{eq:system_WL}
\end{equation}
and $\Cz_{\text{WL}} = E\{\qz_{\text{WL}}(n)\qz_{\text{WL}}^{H}(n)\}$, or 
\begin{align}
    \Cz_{\text{WL}} & =\left[ \begin{array}{cccc}
				         \Cz_{\qz} & \Pz_{\qz} & \Pz_{\qz}^i & \Pz_{\qz}^j\\
					 \Pz_{\qz}^H & \tilde{\Cz}_{\qz} & \tilde{\Cz}_{\qz\qz^i} & \tilde{\Cz}_{\qz\qz^j}\\
					 \Pz_{\qz}^{i H} &  \tilde{\Cz}_{\qz\qz^i}^H & \tilde{\Cz}_{\qz^i} & \tilde{\Cz}_{\qz^i\qz^j}\\
					 \Pz_{\qz}^{j H} &  \tilde{\Cz}_{\qz\qz^j}^H & \tilde{\Cz}_{\qz^i\qz^j}^H & \tilde{\Cz}_{\qz^j}
				      \end{array} \right] 
	\label{eq:correl_WL},			      
\end{align}
where
%
\begin{equation*}
  \begin{array}{ccc}
       \Cz_{\boldsymbol{\alpha}} = E \{\boldsymbol{\alpha} \boldsymbol{\alpha}^H\}, &  \tilde{\Cz}_{\boldsymbol{\alpha}} = E \{\boldsymbol{\alpha}^* \boldsymbol{\alpha}^T\}, & \tilde{\Cz}_{\boldsymbol{\alpha} \boldsymbol{\beta}} = E \{\boldsymbol{\alpha}^* \boldsymbol{\beta}^T\}, \\
       \Pz_{\boldsymbol{\alpha}} = E \{\boldsymbol{\alpha} \boldsymbol{\alpha}^T\}, & \Pz_{\boldsymbol{\alpha}}^i = E \{\boldsymbol{\alpha} \boldsymbol{\alpha}^{iT}\}, & \Pz_{\boldsymbol{\alpha}}^j = E \{\boldsymbol{\alpha} \boldsymbol{\alpha}^{jT}\},\\      
  \end{array}
\end{equation*}
for $\boldsymbol{\alpha}, \boldsymbol{\beta} \in \{\qz(n), \qz^i(n), \qz^j(n)\}$. $\pz_{\text{WL}}$ is given by
\begin{equation*}
  \pz_{\text{WL}} = \text{col}(\hspace{-0.1 cm}\begin{array}{cccc}\pz_\qz, & \pz_{\qz^*}, & \pz_{\qz^{i*}}, & \pz_{\qz^{j*}}\end{array}\hspace{-0.1 cm})
\end{equation*}
and $\pz_{\boldsymbol{\alpha}} = E\{\boldsymbol{\alpha} d^*(n)\}$, $\boldsymbol{\alpha} \in \{\qz(n), \qz^*(n), \qz^{i*}(n), \qz^{j*}(n)\}$.

$\Cz_{\text{WL}}$ and $\pz_{\text{WL}}$ have four times the dimension of their SL counterparts, which are also included in \eqref{eq:system_WL}. 
Using \eqref{eq:system_WL}, we are able to exploit full second-order statistics of the input data, which may not be fully available in \eqref{eq:system_SL}.
Based on statistics provided by $\Cz_{\text{WL}}$ and $\pz_{\text{WL}}$ -- and similarly to the complex case \cite{picimbono1995, tulay2009} -- quaternion 
second-order circularity (also called $\nq$-properness or $\nq$-circularity \cite{santamaria2010}) and the joint quaternion properness can be 
defined as follows.

\subsubsection{$\nq$-properness}

 According to \cite{santamaria2010}, $\qz(n)$ is $\nq$-proper if 
\begin{equation*}
  \tilde{\mathbf{C}}_{\mathbf{q}\mathbf{q}^{i}} = \tilde{\mathbf{C}}_{\mathbf{q}\mathbf{q}^{j}}=\tilde{\mathbf{C}}_{\mathbf{q}^{i}\mathbf{q}^{j}}=\boldsymbol{0}.
\end{equation*}
Noting that
\begin{equation}
    \qz_{\text{WL}}(n)=\Tz \hspace{0.05 cm} \text{col}(\qz(n), \hspace{0.1 cm}  \qz^i(n), \hspace{0.1 cm} \qz^j(n), \hspace{0.1 cm} \qz^k(n)),
    \label{eq:q_wl}
\end{equation}
where
\begin{equation}
      \Tz = \dfrac{1}{2}\left[ \begin{array}{cccc}
				         2 & 0 & 0 & 0\\
					 -1 & 1 & 1 & 1\\
					 1 & -1 & 1 & 1\\
					 1 & 1 & -1 & 1
				      \end{array} \right] \otimes \Iz_N,
	\label{eq:T}		      
\end{equation}
eq. \eqref{eq:correl_WL} can be expressed as
\begin{equation}
    \Cz_{\text{WL}}=   \Tz \left[ \begin{array}{cccc}
				         \Cz_{\qz} & \0z & \0z & \0z\\
					 \0z & \Cz_{\qz}^i  & \0z & \0z\\
					 \0z &  \0z & \Cz_{\qz}^j & \0z\\
					 \0z &  \0z & \0z & \Cz_{\qz}^k
				      \end{array} \right] \Tz^T
	\label{eq:correl_WL2}.	      
\end{equation}
Recalling that
\begin{equation}
\begin{aligned}
	\Cz_{\qz}^\alpha = &E\{\qz^\alpha \qz^{\alpha H}\}=E\{-\alpha \qz \alpha (-\alpha) \qz^H \alpha \}\\
	 = &-\alpha E\{\qz \qz^{H} \}\alpha = -\alpha \Cz_\mathbf{q}\alpha,
    \label{eq:pre_pos}
\end{aligned}
\end{equation}
for $\alpha \in \{ i, j, k\}$, eq. \eqref{eq:correl_WL2} reduces to
\begin{equation}
    \Cz_{\text{WL}} = \tilde{\Tz} (\Iz_4 \otimes \Cz_{\mathbf{q}}) \tilde{\Tz}^H,
    \label{eq:correl_WLn}
\end{equation}
where $ \tilde{\Tz} = \Tz \; \text{diag}(\begin{array}{cccc} 1, & -i, & -j, & -k\end{array}) \otimes \Iz_N$.

From eq. \eqref{eq:correl_WLn}, we notice that for a $\nq$-proper input, only the SL correlation matrix $\Cz_{\qz}$ is required
to access the second-order statistics of $\qz(n)$, similarly to the complex-proper case reported in \cite{picimbono1995}.

\subsubsection{Joint $\nq$-properness}
$\nq$-properness is related only to the input signal $\qz(n)$. 
However, $d(n)$ and $\qz(n)$ can also be jointly $\nq$-proper \cite{santamaria2010}, 
which implies the additional restrictions
\begin{equation}
  E\{\qz(n) d^{\alpha*}(n)\}=\0z, \; \forall \alpha \in \{ i,\; j, \; k\}.
  \label{eq:cond}
\end{equation}
In this case, SL and WL estimation provide the same result.
To show this equivalence, use \eqref{eq:q_wl} to write $\pz_{\text{WL}}$ as given by
\begin{equation}
    \pz_{\text{WL}} = \Tz \left[ \begin{array}{c} \pz_\qz \\ -i E\{\qz(n) d^{i*}(n)\} i \\ -j E\{\qz(n) d^{j*}(n)\} j  \\ -k E\{\qz(n) d^{k*}(n)\} k  \end{array}\right].
    \label{eq:p_longo}
\end{equation}
Using eq. \eqref{eq:cond} in \eqref{eq:p_longo}, we obtain
\begin{equation*}
    \pz_{\text{WL}} = \Tz \hspace{0.1 cm} \text{col}(\pz_\qz, \hspace{0.1 cm} \0z, \hspace{0.1 cm} \0z, \hspace{0.1 cm} \0z).
\end{equation*}
Finally, use the joint $\nq$-properness restrictions to express \eqref{eq:system_WL} as 
\begin{equation}
       \Tz \hspace{-0.1 cm}\left[ \begin{array}{cccc}
				         \Cz_{\qz} & \0z & \0z & \0z\\
					 \0z & \Cz_{\qz}^i  & \0z & \0z\\
					 \0z &  \0z & \Cz_{\qz}^j & \0z\\
					 \0z &  \0z & \0z & \Cz_{\qz}^k
				      \end{array} \right] \hspace{-0.1 cm} \Tz^T  \wz_{\text{WL}}(n) = \Tz \hspace{-0.1 cm}\left[ \begin{array}{c} \pz_\qz \\ \0z \\ \0z \\  \0z \end{array} \right]
	\label{eq:correl_WL3}.		      
\end{equation}
To solve \eqref{eq:correl_WL3}, define 
\begin{equation*}
  \wz_{\text{WL}}=\text{col}(\hspace{-0.1 cm}\begin{array}{cccc} \wz_{\text{WL}, 1}, & \wz_{\text{WL}, 2} & \wz_{\text{WL}, 3} & \wz_{\text{WL}, 4}\end{array}\hspace{-0.1 cm}),
\end{equation*}
where we drop the time indices to simplify notation. Noting that $\Tz$ is invertible, and left-multiplying \eqref{eq:correl_WL3} by $\Tz^{-1}$, we obtain the following systems of equations
\begin{align}
    &\Cz_\qz ( 2\wz_{\text{WL}, 1} - \wz_{\text{WL}, 2} + \wz_{\text{WL}, 3}+ \wz_{\text{WL}, 4}) = 2 \pz_q,
    \label{eq1}\\
   &\Cz_\qz^i (\wz_{\text{WL}, 2} - \wz_{\text{WL}, 3} + \wz_{\text{WL}, 4} ) = \0z,   	
  \label{eq2}\\
    &\Cz_\qz^j (\wz_{\text{WL}, 2} + \wz_{\text{WL}, 3} - \wz_{\text{WL}, 4} ) = \0z,
    \label{eq3}\\
    &\Cz_\qz^k (\wz_{\text{WL}, 2} + \wz_{\text{WL}, 3} + \wz_{\text{WL}, 4} ) = \0z.
    \label{eq4}
\end{align}
The only solution to \eqref{eq2}, \eqref{eq3} and \eqref{eq4} is trivial, which is easy to verify. For this purpose, consider eq. \eqref{eq2} and define 
\begin{equation*}
  \bz = \wz_{\text{WL}, 2} - \wz_{\text{WL}, 3} + \wz_{\text{WL}, 4}.
\end{equation*}
Recall eq. \eqref{eq:pre_pos}, and assume that $\Cz_\qz$ is full-rank (which is normally the case). The null-space 
\cite{meyerMAALA00} of $\Cz_\qz^i$ has only the trivial solution $\bz=\0z$. Applying the same argument to \eqref{eq3} and \eqref{eq4}, the trivial 
solution is straightforward. Using these results, we obtain a system of equations to compute $\wz_{\text{WL}, 2}$, $\wz_{\text{WL}, 3}$ and $\wz_{\text{WL}, 4}$, i.e.,
\begin{equation*}
  \left( \left[ \begin{array}{ccc}
      1 & -1 & 1\\
      1 & 1 & -1\\
      1 & 1 & 1
  \end{array}\right] \otimes \Iz_N\right)   \left[ \begin{array}{c}
      \wz_{\text{WL}, 2}\\
      \wz_{\text{WL}, 3}\\
      \wz_{\text{WL}, 4}
  \end{array}\right]  =  \left[ \begin{array}{c}
      \0z\\
      \0z\\
      \0z
  \end{array}\right],
\end{equation*}
whose only solution is the null vector. Substituting $\wz_{\text{WL}, 2}=\wz_{\text{WL}, 3}=\wz_{\text{WL}, 4}=\0z$ in 
\eqref{eq1}, only the system of equations of \eqref{eq:system_SL} must be computed. This result shows that the SL and WL approaches 
are equivalent in terms of MSE performance when the input and the desired signal are jointly $\nq$-proper.

In summary, when $\qz(n)$ and $d(n)$ are joint $\nq$-proper, a WL estimator is in general not advantageous, since the MSE performance is the same of a SL 
approach, but the number of computations required to obtain the solution is higher. On the other hand, when joint $\nq$-properness does not hold, 
a WL-based technique can exploit full second-order statistics of the input data to solve \eqref{eq:system_WL}, improving the MSE performance.



\section{Quaternion gradients}
\label{sec:quaternion_gradients}


In this section, we propose a general definition for the quaternion derivatives, which is used to study the convergence of WL-QLMS-based 
algorithms and to obtain the fastest-converging WL-QLMS-like algorithm in two situations: 
1) When at most two of the quaternion elements in the regressor are correlated; and
2) When the WL algorithms use a real regressor vector, as proposed in \cite{fernandoSSP2011}. 
We also use the general description for the gradient to develop mean and mean-square analyses of real-regressor vector WL-QLMS algorithms, 
providing accurate tools to design and evaluate the performance.

To obtain a general quaternion gradient, start with the definition of the cost function
\begin{equation*}
    f(\wz)=e(n)e^*(n),
\end{equation*}
where $e(n)$ corresponds to the error between a desired quantity $d(n)$ and the estimated value
\begin{equation*}
   \hat{d}(n)=\wz^H(n)\qz(n),
\end{equation*}
and $\wz = \wz_\text{R} + i\wz_\text{I} + j\wz_\text{J} + k\wz_\text{K}$.
Regardless the method to compute the gradient, and based on the isomorphism between $\nq$ and $\nr^4$, the quaternion gradients proposed in the literature have the general form
\begin{equation}
    \nabla_{\wz}f = a\dfrac{\partial f}{\partial \wz_\text{R}}+ i b\dfrac{\partial f}{\partial \wz_\text{I}} + j c\dfrac{\partial f}{\partial \wz_\text{J}} +  k  d\dfrac{\partial f}{\partial \wz_\text{K}},
    \label{eq:deriv}
\end{equation}
for real $\{a,b,c,d\}$, where
\begin{equation}
  \dfrac{\partial f}{\partial \wz_\alpha}= \dfrac{\partial e(n)}{\partial \wz_\alpha}e^*(n) + e(n)\dfrac{\partial e^*(n)}{\partial \wz_\alpha}, 
  \label{eq:partial_eq}
\end{equation}
for $\alpha \in \{R,I,J,K\}$. Using eq. \eqref{eq:partial_eq} in \eqref{eq:deriv}, and defining
\begin{equation}
    \begin{aligned}
    g = (a + b + c + d)\\
    h = (a - b - c - d)
    \end{aligned}
  \label{eq:h_e_g}
\end{equation}
eq. \eqref{eq:deriv} becomes
\begin{equation}
    \nabla_{\wz}f = -g \qz(n)e^*(n)  - h e(n)\qz^*(n).
    \label{eq:grad}
\end{equation}
From \eqref{eq:grad}, one can check that the quaternion gradient of \cite{Took2011} is obtained when $h=1/2$ ad $g = -1/4$, 
and that the i-gradient of \cite{mandic_iQLMS} appears when $h=3/4$ ad $g = 0$. We also note that the reviewed quaternion gradient, 
presented in \cite{iqlms_review}, is obtained when $h=2$ and $g=0$.
Considering eq. \eqref{eq:h_e_g}, it is easy to note that the same $h$ and $g$ can be obtained for different values of $a$, $b$, $c$ and $d$. 
Moreover, note that gradient-based algorithms have a general update law given by \cite{sayed2008adaptive}
\begin{equation}
  \wz(n+1) = \wz(n) - \mu \nabla_{\wz} f,
  \label{eq:update}
\end{equation}
where $\mu$ is the step-size. In this case, from the point of view of the algorithm, it is not important if the gradient uses $h$ and $g$ or scaled versions of them -- as long as they are both scaled by the same value -- since the scale factor can be absorbed by $\mu$. This fact emphasizes that different gradients and quaternion derivatives can be applied to define the same QLMS-like algorithm.

Using eq. \eqref{eq:grad} in \eqref{eq:update}, we define the general form of a QLMS-like algorithm, that is
\begin{equation}
    \wz(n+1) = \wz(n) + \mu \left[ g \qz(n) e^*(n) + h e(n) \qz^*(n) \right].
    \label{eq:qlms_gen}
\end{equation}
Noting that
\begin{multline}
    g \qz(n) e^*(n) + h e(n) \qz^*(n) = \\=(g + h)\qz(n) e^*(n) - 2h \imag\{\qz(n) e^*(n)\},
    \label{eq:relation}
\end{multline}
and  using eq. \eqref{eq:relation} in \eqref{eq:qlms_gen}, one gets
\begin{multline}
    \wz(n+1) = \wz(n) + \mu \left[ (g + h)\qz(n) e^*(n)\right. \\ \left. - 2h \Im\{\qz(n) e^*(n)\} \right].
    \label{eq:qlms_gen1}
\end{multline}
In order to study the mean behavior of \eqref{eq:qlms_gen1}, define the optimum set of
coefficients $\wz_o$, and assume that $d(n)$ can be modeled as given by
\begin{equation}
  d(n) = \wz_o^H \qz(n) + v(n),
  \label{eq:desejado}
\end{equation}
where the elements of $v(n)$ are i.i.d. Gaussian noise, independent of $\qz(n)$, $E\{v(n)\}=0$, and $E\{|v(n)|^2\}=\sigma_v^2$.
Subtracting  \eqref{eq:qlms_gen1} from $\wz_o$ and taking the expectation, one gets
\begin{multline}
    \bar{\wz}(n+1)=\bar{\wz}(n)-\mu E\{(g+h)\qz(n)e^*(n)  \\ + 2h \Im\{\qz(n) e^*(n)\} \},
    \label{eq:qlms_gen2}
\end{multline}
where  $\bar{\wz}(n) = E\{\tilde{\wz}(n)\}$ and
\begin{equation}
    \tilde{\wz}(n) = \wz_o - \wz(n).
    \label{eq:w_tilde}
\end{equation}
Use \eqref{eq:desejado} in \eqref{eq:qlms_gen1}, and assume the independence approximation usually applied to study LMS-like 
algorithms \cite{sayed2008adaptive}
\begin{equation*}
   E\{\qz(n)\qz^H(n)\tilde{\wz}(n)\} \approx E\{\qz(n)\qz^H(n)\}\bar{\wz}(n).
\end{equation*}
After some algebra manipulation, eq. \eqref{eq:qlms_gen2} becomes
\begin{equation}
\begin{aligned}
    &\bar{\wz}(n+1) = \bar{\wz}(n)  
                     \\ &- \mu (g + h) \left[ E\{\qz(n)\qz^H(n)\}\bar{\wz}(n) + E\{\qz(n)v^*(n)\} \right] 
	             \\ &+ 2\mu h \imag \left\lbrace  E\{\qz(n)\qz^H(n)\}\bar{\wz}(n) +  E\{\qz(n)v^*(n)\}\right\rbrace.  
    \label{eq:w_quat}
\end{aligned}
\end{equation}
Considering the assumptions for $v(n)$, eq. \eqref{eq:w_quat} reduces to
\begin{multline}
    \bar{\wz}(n+1) = \bar{\wz}(n)  
                      - \mu (g + h) E\{\qz(n)\qz^H(n)\}\bar{\wz}(n)
	             \\ + 2\mu h \imag\left\lbrace  E\{\qz(n)\qz^H(n)\}\bar{\wz}(n)\right\rbrace.  
    \label{eq:equation_quat}
\end{multline}
With no loss of generality, consider from this point on that we are performing the analysis of a WL algorithm. In this case, we add a subscript WL 
to the variables to write
\begin{multline}
    \bar{\wz}_{\text{WL}}(n+1) = \bar{\wz}_{\text{WL}}(n)  \\
                      - \mu_{\text{WL}} (g + h) E\{\qz_{\text{WL}}(n)\qz_{\text{WL}}^H(n)\}\bar{\wz}_{\text{WL}}(n)
	             \\ + 2\mu_{\text{WL}} h \imag \left\lbrace  E\{\qz_{\text{WL}}(n)\qz_{\text{WL}}^H(n)\}\bar{\wz}_{\text{WL}}(n)\right\rbrace. 
    \label{eq:mean_wl}
\end{multline}
In terms of $\mathbf{\bar{w}}_{\text{WL}}(n)$, this is a nonlinear recursion (due to the $\imag \{ \cdot \}$ operator). 
In order to obtain a linear recursion, we resort to the extended variables.
Using the extended entities (as proposed in \cite{fernandoSSP2011}), define the extended $\bar{\wz}_{\text{WL}}(n)$, i.e.,
\begin{equation*}
    \bar{\wz}_{\text{ext}}(n)= \text{col}(\bar{\wz}_{{\text{WL}}_{\text{R}}}(n), \hspace{0.1 cm} \bar{\wz}_{{\text{WL}}_{\text{I}}}(n),  \hspace{0.1 cm} \bar{\wz}_{{\text{WL}}_{\text{J}}}(n), \hspace{0.1 cm} \bar{\wz}_{{\text{WL}}_{\text{K}}}(n))
\end{equation*}
and the extended version of $E\{\qz_{\text{WL}}(n)\qz_{\text{WL}}^H(n)\}\bar{\wz}_{\text{WL}}(n)$, which is given by
$\Cz_{\text{ext}}\bar{\wz}_{\text{ext}}(n)$,
with
\begin{equation}
    \Cz_{\text{ext}}=\left[\begin{array}{cccc} 
	 \RM & \IM^T & \JM^T & \KM^T\\
	  \IM &\RM & \KM^T & \JM\\
	  \JM & \KM &\RM & \IM^T\\
	  \KM & \JM^T & \IM &\RM
    \end{array} \right].
    \label{eq:real_c}
\end{equation}
The matrix $\RM$ contains the real elements of $E\{\qz_{\text{WL}}(n)\qz_{\text{WL}}^H(n)\}$, and $\IM$, $\JM$ and $\KM$ 
are the imaginary parts for $i$, $j$ and $k$, respectively. $\RM$ is a symmetric matrix, and 
$\IM=-\IM^T$, $\JM=-\JM^T$ and $\KM=-\KM^T$. Defining the WL regressor vector as
\begin{equation}
	\mathbf{q}_{\text{WL}}(n) = \mathbf{q}_{{\text{WL}}_{\text{R}}}(n) + i\mathbf{q}_{{\text{WL}}_{\text{I}}}(n) + j\mathbf{q}_{{\text{WL}}_{\text{J}}}(n) + k\mathbf{q}_{{\text{WL}}_{\text{K}}}(n),
\end{equation}
these matrices are given by
\begin{align*}
	\RM & = \mathbf{C}_\text{R} + \mathbf{C}_{\text{I}} + \mathbf{C}_\text{J} + \mathbf{C}_\text{K} \\
	\IM & = -\mathbf{C}_{\text{RI}} + \mathbf{C}^T_{\text{RI}} - \mathbf{C}_{\text{JK}} + \mathbf{C}^T_{\text{JK}}, \\
	\JM & = -\mathbf{C}_{\text{RJ}} + \mathbf{C}^T_{\text{RJ}} + \mathbf{C}_{\text{IK}} -\mathbf{C}^T_{\text{IK}},\\
	\KM & = -\mathbf{C}_{\text{RK}} + \mathbf{C}^T_{\text{RK}} - \mathbf{C}_{\text{IJ}} +\mathbf{C}^T_{\text{IJ}},
\end{align*}
where
\begin{align}
&\mathbf{C}_{\text{R}} = E\{\mathbf{q}_{{\text{WL}}_{\text{R}}}(n) \mathbf{q}_{{\text{WL}}_{\text{R}}}^T(n)\}, && \hspace{-0.1 cm}\mathbf{C}_{\text{I}} = E\{\mathbf{q}_{{\text{WL}}_{\text{I}}}(n) \mathbf{q}_{{\text{WL}}_{\text{I}}}^T(n)\}, \nonumber \\ 
&\mathbf{C}_{\text{J}} = E\{\mathbf{q}_{{\text{WL}}_{\text{J}}}(n) \mathbf{q}_{{\text{WL}}_{\text{J}}}^T(n)\}, && \hspace{-0.1 cm}\mathbf{C}_{\text{K}} = E\{\mathbf{q}_{{\text{WL}}_{\text{K}}}(n) \mathbf{q}_{{\text{WL}}_{\text{K}}}^T(n)\}, \nonumber \\ 
&\mathbf{C}_{\text{RI}} = E\{\mathbf{q}_{{\text{WL}}_{\text{R}}}(n) \mathbf{q}_{{\text{WL}}_{\text{I}}}^T(n)\}, && \hspace{-0.1 cm}\mathbf{C}_{\text{JK}} = E\{\mathbf{q}_{{\text{WL}}_{\text{J}}}(n) \mathbf{q}_{{\text{WL}}_{\text{K}}}^T(n)\}, \nonumber \\ 
&\mathbf{C}_{\text{RJ}} = E\{\mathbf{q}_{{\text{WL}}_{\text{R}}}(n) \mathbf{q}_{{\text{WL}}_{\text{J}}}^T(n)\}, && \hspace{-0.1 cm}\mathbf{C}_{\text{IK}} = E\{\mathbf{q}_{{\text{WL}}_{\text{I}}}(n) \mathbf{q}_{{\text{WL}}_{\text{K}}}^T(n)\}, \nonumber \\ 
&\mathbf{C}_{\text{RK}} = E\{\mathbf{q}_{{\text{WL}}_{\text{R}}}(n) \mathbf{q}_{{\text{WL}}_{\text{K}}}^T(n)\}, && \hspace{-0.1 cm}\mathbf{C}_{\text{IJ}} = E\{\mathbf{q}_{{\text{WL}}_{\text{I}}}(n) \mathbf{q}_{{\text{WL}}_{\text{J}}}^T(n)\}. \nonumber
\end{align}
Using the extended entities, $\imag \{ E\{\qz_{\text{WL}}(n)\qz_{\text{WL}}^H(n)\}\bar{\wz}_{\text{WL}}(n) \}$
is replaced by $\Cz^{\text{Im}}_{\text{ext}}\bar{\wz}_{\text{ext}}(n)$, where,
\begin{equation}
    \Cz^{\text{Im}}_{\text{ext}}= \left[ \text{diag}(0,\; 1,\; 1,\; 1)\otimes \mathbf{I}_{4N} \right] \Cz_{\text{ext}}.
    \label{eq:imag_c}
\end{equation}
Applying eqs. \eqref{eq:real_c} and \eqref{eq:imag_c}, the extended version of \eqref{eq:mean_wl} is given by
\begin{multline*}
    \bar{\wz}_{\text{ext}}(n+1) = \bar{\wz}_{\text{ext}}(n) 
                      - \mu_{\text{WL}} (g + h) \Cz_{\text{ext}}\bar{\wz}_{\text{ext}}(n)
	             \\ + 2\mu_{\text{WL}} h \Cz^{\text{Im}}_{\text{ext}}\bar{\wz}_{\text{ext}}(n).
\end{multline*}
Considering the structure of $\Cz_{\text{ext}}$ and $\Cz^{\text{Im}}_{\text{ext}}$, one can write
\begin{equation}
    \bar{\wz}_{\text{ext}}(n+1) = \left( \Iz_{4N} - \mu_{WL} \Gz_{\text{ext}} \Cz_{\text{ext}}\right) \bar{\wz}_{\text{ext}}(n),
    \label{eq:mean_wl2}
\end{equation}
and 
\begin{equation}
      \Gz_{\text{ext}} = \Gz \otimes \Iz_{4N}, 
      \label{eq:G_ext}
\end{equation}
\begin{equation*}
	\Gz = \text{diag}(\hspace{-0.1 cm}\begin{array}{c} 
			(g+h), (g-h), (g-h), (g-h)
			\end{array}\hspace{-0.1 cm}).
\end{equation*}
The product $\Gz_{\text{ext}} \Cz_{\text{ext}}$ is the extended quaternion correlation matrix, which is obtained when the
quaternion entries of \eqref{eq:mean_wl} are mapped to the real field. $\Gz_{\text{ext}} \Cz_{\text{ext}}$ can be used to study the convergence of the 
algorithms, which is assessed through its eigenvalue spread. In Appendix \ref{ap:GC}, we show that $\Cz_{\text{ext}}$ is a positive 
semi-definite matrix \cite{hornTMA91}, and that $\Gz_{\text{ext}}$ must be at least positive semi-definite to avoid divergence in 
\eqref{eq:mean_wl2}.
In the next section, the convergence of \eqref{eq:mean_wl2} is studied in two situations, where we prove that a gradient with $g=0$ 
provides the faster-converging algorithms. We start with the case when at most two quaternion elements are correlated, and then we study 
quaternion algorithms using a real regressor vector, as proposed in \cite{fernandoSSP2011}.

\subsection{Case one: signals with correlation between at most two quaternion elements}
\label{sec:two_quaternions}

Define $\alpha = (g + h)/(g-h)$ and assume that $\JM = \KM =\boldsymbol{0}$, such that
$\mathbf{G}_\text{ext} \mathbf{C}_\text{ext}$ is written as
\begin{equation}
	(g-h) \left[ \begin{array}{cccc}
		\mathbf{A}(\alpha) & \boldsymbol{0}\\
		\boldsymbol{0} & \mathbf{B}
	\end{array}\right] = (g-h) \mathbf{C}(\alpha), 
	\label{eq:case_3}
\end{equation}
where
\begin{equation*}
	\mathbf{A}(\alpha)= \left[ \begin{array}{cc}
		\alpha\RM & \alpha\IM^T\\
		\IM & \RM
	\end{array}\right]
\end{equation*}
and $\mathbf{B}=\mathbf{A}(1)$. In addition, define the spreading factor of $\mathbf{A}(\alpha)$, $\mathbf{B}$ and 
$\mathbf{C}(\alpha)$ as $SF_{\mathbf{A}}$, $SF_{\mathbf{B}}$ and $SF_{\mathbf{C}}$, respectively. We want to show that the smallest condition number 
of \eqref{eq:case_3} is obtained when $\alpha=1$. For this purpose, we first consider $\mathbf{A}(\alpha)$ to show that it has the minimum 
condition number when $\alpha=1$. Then, we show that the smallest condition number of $\mathbf{C}(\alpha)$ is also obtained when $\alpha=1$.

First, consider $\mathbf{A}(\alpha)$. One can show that $\mathbf{A}(\alpha)$ corresponds to a row scaling of $\mathbf{B}$, since it can be written as 
\begin{equation*}
\mathbf{A}(\alpha)= \mathbf{D}(\alpha)\mathbf{B},
\end{equation*}
where $\mathbf{D}(\alpha)$ is the scaling matrix, given by $\mathbf{D}(\alpha) = \text{diag}(\alpha, \; 1) \otimes \mathbf{I}_{8N}$.
In \cite{braatz94}, the problem of determining the row scaling of a matrix which minimizes the Euclidean condition number is addressed. 
It is shown that this problem is convex, and can be solved by convex optimization. With this information, we know that there is a value of 
$\alpha$ which minimizes the condition number, which we now find out. 

Studying matrix $\mathbf{A}(\alpha)$, one can show that the condition number of $\mathbf{A}(\alpha)$ is equal to that of $\mathbf{A}(1/\alpha)$. 
For this purpose, we use an unitary transformation \cite{meyerMAALA00} -- which does not change the eigenvalues of a matrix -- to show that the 
spreading factors of $\mathbf{A}(\alpha)$ and $\mathbf{A}(1/\alpha)$ have the same value, since
\begin{equation*}
	\left[ \begin{array}{cc} 
			\boldsymbol{0} & \mathbf{I}_{8N}\\
			-\mathbf{I}_{8N} & \boldsymbol{0}	
	\end{array}\right] 
	\mathbf{A}(\alpha) 
	\left[ \begin{array}{cc} \boldsymbol{0} & -\mathbf{I}_{8N}\\
							 \mathbf{I}_{8N} & \boldsymbol{0}	
	\end{array}\right] = \alpha \mathbf{A}(1/ \alpha), 
\end{equation*}
where we also have used the fact that $\IM=-\IM^T$. Notice that when $\alpha \rightarrow \infty$ or when $\alpha \rightarrow 0$, 
the condition numbers of $\mathbf{A}(\alpha)$ and $\mathbf{A}(1/\alpha)$ both go to $\infty$. 
Since $SF_{\mathbf{A}}(\alpha) = SF_{\mathbf{A}}(1/\alpha)$, and $SF_{\mathbf{A}}(\alpha)$ is a convex function \cite{braatz94}, 
the minimum condition number must be at the point which defines the axis of symmetry of the problem. This point corresponds to $\alpha=1$, 
where $SF_{\mathbf{A}}= SF_{\mathbf{B}}$. For different values of $\alpha$, $SF_{\mathbf{A}} > SF_{\mathbf{B}}$.

Using the result for $\mathbf{A}(\alpha)$, we can evaluate the condition number of matrix $\mathbf{C}$ in eq. \eqref{eq:case_3}. 
 $SF_\mathbf{C}$ is given by\footnote{Note that the constant $(g-h)$ does not affect the eigenvalue spread of eq. \eqref{eq:case_3}, 
since it multiplies all the eigenvalues.} 

\begin{equation}
	SF_{\mathbf{C}} = \frac{\text{max}\left\lbrace \lambda_{\text{Max}}(\mathbf{A}(\alpha)),  \lambda_{\text{Max}}(\mathbf{B})\right\rbrace }{\text{min}\left\lbrace \lambda_{\text{Min}}(\mathbf{A}(\alpha)),  \lambda_{\text{Min}}(\mathbf{B})\right\rbrace},
	\label{eq:SFC}
\end{equation}  
where $\lambda_{\text{Max}}(\mathbf{A}(\alpha))$ and  $\lambda_{\text{Max}}(\mathbf{B})$ are the maximum eigenvalues of 
$\mathbf{A}(\alpha)$ and $\mathbf{B}$, respectively, and $\lambda_{\text{Min}}(\mathbf{A}(\alpha))$ and  $\lambda_{\text{Min}}(\mathbf{B})$ are
 the minimum eigenvalues. Using \eqref{eq:SFC}, it is possible to list the cases which can appear in the computation of $SF_{\mathbf{C}}$, i.e.,
\begin{enumerate}
	\item $\lambda_{\text{Max}}(\mathbf{A}(\alpha)) > \lambda_{\text{Max}}(\mathbf{B})$ and $\lambda_{\text{Min}}(\mathbf{A}(\alpha)) \leq \lambda_{\text{Min}}(\mathbf{B})$. In this case, $SF_{\mathbf{C}}= \lambda_{\text{Max}}(\mathbf{A}(\alpha))/\lambda_{\text{Min}}(\mathbf{A}(\alpha)) > SF_{\mathbf{B}}$.
	\item $\lambda_{\text{Max}}(\mathbf{A}(\alpha)) > \lambda_{\text{Max}}(\mathbf{B})$ and $\lambda_{\text{Min}}(\mathbf{A}(\alpha)) \geq \lambda_{\text{Min}}(\mathbf{B})$. In this case, $SF_{\mathbf{C}}= \lambda_{\text{Max}}(\mathbf{A}(\alpha))/\lambda_{\text{Min}}(\mathbf{B}) > SF_{\mathbf{B}}$.
	\item  $\lambda_{\text{Max}}(\mathbf{A}(\alpha)) \leq \lambda_{\text{Max}}(\mathbf{B})$ and $\lambda_{\text{Min}}(\mathbf{A}(\alpha)) <  \lambda_{\text{Min}}(\mathbf{B})$. In this case, $SF_{\mathbf{C}}= \lambda_{\text{Max}}(\mathbf{B})/\lambda_{\text{Min}}(\mathbf{A}(\alpha)) > SF_{\mathbf{B}}$.
	\item $\lambda_{\text{Max}}(\mathbf{A}(\alpha)) = \lambda_{\text{Max}}(\mathbf{B})$ and $\lambda_{\text{Min}}(\mathbf{A}(\alpha)) = \lambda_{\text{Min}}(\mathbf{B})$. In this case, $SF_{\mathbf{C}}=  SF_{\mathbf{B}}$. 
\end{enumerate}
Notice that the condition $\lambda_{\text{Max}}(\mathbf{A}(\alpha)) < \lambda_{\text{Max}}(\mathbf{B})$ and $\lambda_{\text{Min}}(\mathbf{A}(\alpha))
 > \lambda_{\text{Min}}(\mathbf{B})$ is not possible, since we have shown that $SF_{\mathbf{A}} \geq SF_{B}$ always.
 
For the three first conditions, the eigenvalue spread is always increased, while the last possibility reveals that the minimum value of $SF_\mathbf{C}$ 
occurs when $\mathbf{A}(\alpha)=\mathbf{B}$, for $\alpha = 1$. Recalling that $\alpha=(g+h)/(g-h)$, we conclude that $h$ must be zero, 
showing that the minimum eigenvalue spread is obtained with a gradient similar to that of \cite{mandic_iQLMS} or \cite{iqlms_review}. 

Note that a similar approach can be used when only $\RM$ and $\JM$ or $\RM$ and $\KM$ are different from a zero matrix, showing that 
the minimum eigenvalue spread is also obtained with $\alpha = 1$.

\subsubsection{Particular case: all the quaternion elements are uncorrelated among them}
\label{sec:uncor_input}

Consider now the particular case when all the quaternion elements are uncorrelated, which means that $\IM=\JM=\KM = \0z$. In this situation,
$\Cz_\text{ext}=\Iz_4 \otimes \RM$, and one can easily show that 
$\mathbf{A}(\alpha) = \text{diag}(\begin{array}{cc} \alpha, & 1 \end{array}) \otimes \RM$ and $\mathbf{B}=\mathbf{A}(1)$,
such that the result of Section \ref{sec:two_quaternions} still holds.   

%

\subsection{Case two: Widely-linear algorithms using real data vector}
\label{sec:real_regressor}

Consider now a WL quaternion algorithm implemented with a real-regressor vector 
\begin{equation}
    \xz(n) = \text{col}(\qz_R(n), \hspace{0.1 cm} \qz_I(n), \hspace{0.1 cm} \qz_J(n), \hspace{0.1 cm} \qz_k(n)),
    \label{eq:real_regressor}
\end{equation}
as proposed for RC-WL-QLMS in \cite{fernandoSSP2011}. For this aproach, the correlation matrix $\Cz_{\xz}$ is real and given by
\begin{equation}
  \Cz_{\xz} = E\{\xz(n)\xz^T(n)\},
  \label{eq:Cx}
\end{equation}
and \eqref{eq:mean_wl2} changes to
\begin{equation}
    \bar{\wz}_{\text{ext}}(n+1) = \left( \Iz_{4N} - \mu_{WL} \Gz_{\text{ext}} \Cz_{\xz_{\text{ext}}}\right) \bar{\wz}_{\text{ext}}(n),
    \label{eq:mean_wl3}
\end{equation}
where
\begin{equation}
    \Cz_{\xz_{\text{ext}}}=\Iz_4 \otimes \Cz_{\xz}.  
    \label{eq:Cx_ext}
\end{equation}
The convergence is studied through the product
\begin{equation*}
    \Cz_{\xz_{\text{ext}}}\Gz_{\text{ext}} = (\Iz_4 \otimes \Cz_{\xz})(\Gz \otimes \Iz_{4N})=\Gz \otimes \Cz_{\xz}.
\end{equation*}
Using an approach similar to that of Section \ref{sec:two_quaternions}, we compute
\begin{equation}
    SF_{\Gz_{\text{ext}} \mathbf{C}_{{\mathbf{x}}_\text{ext}}} = \text{max}\left\lbrace \dfrac{(g+h)}{(g-h)}, \dfrac{(g-h)}{(g+h)}\right\rbrace \times \dfrac{\lambda_{\text{Max}}(\mathbf{C}_{\xz})}{\lambda_{\text{Min}}(\mathbf{C}_{\xz})},
    \label{eq:SF}
\end{equation}
where $\lambda_{\text{Max}}(\mathbf{C}_{\xz})$ and $\lambda_{\text{Min}}(\mathbf{C}_{\xz})$ are the largest and the smallest eigenvalues of 
$\Cz_{\xz}$, respectively. The best choice for $g$ and $h$ is again obtained when $(g+h)/(g-h)=1$, which results in $h=0$ and $a = b+c+d$. 
 
Notice that the family of gradients which uses $g=-1/2$ and $h = 1/4$ (or scaled versions of $h$ and $g$) leads to the definition of the RC-WL-QLMS 
algorithm of \cite{fernandoSSP2011}. 
However, these choices for $g$ and $h$ do not lead to the minimum condition number of $\Gz_{\text{ext}} \Cz_{\xz_{\text{ext}}}$. When $h=0$, eq. \eqref{eq:SF} is reduced 
to its minimum value, such that the condition number depends only on the eigenvalue spread of 
$\mathbf{C}_{\xz}$. In this case, any gradient such that $h=0$ will result in an algorithm which achieves the fastest-converging algorithm using $\xz(n)$.
In addition, as presented in \cite{fernandoSSP2011}, the use of real-regressor vectors leads to lower-complexity quaternion algorithms, since many 
operations can be avoided. Based on these ideas, in Section \ref{sec:rc_wl_qlms} we propose the fastest-converging WL-QLMS algorithm with a regressor 
vector given by \eqref{eq:real_regressor}, which uses even less computations than RC-WL-QLMS \cite{fernandoSSP2011}.

We must emphasize an important aspect of the analysis presented here. This approach is valid for both correlated and uncorrelated input, since 
no initial restriction to the correlation matrix is necessary to the study of the eigenvalues of $\mathbf{C}_{{\text{X}_{\text{ext}}}}$. 
Moreover, it can be used to any WL quaternion algorithm with a real-regressor vector obtained with the concatenation of the elements of 
$\mathbf{q}(n)$.

In the following sections, we use the general equation proposed to describe the quaternion gradients to obtain 
a new reduced-complexity algorithm with real regressor vector. This general description also allows us
to develop an analysis for WL-QLMS algorithms with real data vector.

\section{Proposed new reduced-complexity WL-QLMS algorithm with real-regressor vector}
\label{sec:rc_wl_qlms}


In reference \cite{fernandoSSP2011}, the RC-WL-QLMS algorithm was proposed as a lower-complexity alternative to WL-QLMS. In that paper, 
the algorithm was not formally defined in terms of a gradient, but using the general approach of eq. \eqref{eq:qlms_gen}, one can 
show that RC-WL-QLMS corresponds to a real-regressor WL-QLMS algorithm where $h=-1/4$ and $g=1/2$ (or where $h$ and $g$ are multiples of these values). The algorithm is presented in 
Table \ref{tab:rc_wl_qlms}, where $\hat{d}_{\text{RC}_{\text{Q}}}(n)$ is the estimate of $d(n)$, $e_{\text{RC}_{\text{Q}}}(n)$ is the error, 
$\wz_{\text{RC}_{\text{Q}}}(n)$ are 
the weights, and $\mu_{\text{RC}_{\text{Q}}}$ is the step-size. 
\begin{table}[htb]
 \centering
   \caption{\it RC-WL-QLMS algorithm}
   \begin{tabular}{|l|} 
      \hline
      $\hat{d}_{\text{RC}_{\text{Q}}}(n)=\wz_{\text{RC}_{\text{Q}}}^{H}(n)\xz(n)$ \\
      $e_{\text{RC}_{\text{Q}}}(n)= d(n)-\hat{d}_{\text{RC}_{\text{Q}}}(n)$\\
      $\wz_{\text{RC}_{\text{Q}}}(n+1)=\wz_{\text{RC}_{\text{Q}}}(n)+\mu_{\text{RC}_{\text{Q}}} \left[ \frac{e_{\text{RC}_{\text{Q}}}^{*}(n)}{2}  - \frac{e_{\text{RC}_{\text{Q}}}(n)}{4}\right] \xz(n)$ \\
      \hline
   \end{tabular}
   \label{tab:rc_wl_qlms}
\end{table}

It is shown in \cite{fernandoSSP2011} and \cite{fernandoISWCS2012} that the RC-WL-QLMS algorithm converges faster than WL-QLMS, and that it is also less costly to compute. 
Yet, our general approach to the gradient indicates that RC-WL-QLMS is not the fastest-converging WL-QLMS algorithm with real-regressor vector, 
since the gradient used to define it uses $h \neq 0$. We exploit this fact to propose a new reduced-complexity algorithm, which outperforms RC-WL-QLMS 
both in the computational cost and in the convergence rate. For this purpose, assume that the real and imaginary parts of the original SL data vector are 
stacked up in the vector $\mathbf{x}(n)$, as expressed in eq. \eqref{eq:real_regressor}. Define the estimate $\hat{d}_{\text{RC}}(n)$ of the desired
signal $d(n)$, i.e.,
\begin{equation}
    \hat{d}_{\text{RC}}(n)=\wz_{\text{RC}}^H(n)\xz(n),
    \label{eq:drc}
\end{equation}
where $\wz_{\text{RC}}(n)$ is the vector of weights. The error is given by 
\begin{equation}
    e_{\text{RC}}(n) = d(n) - \hat{d}_{\text{RC}}(n).
    \label{eq:erc}
\end{equation}
Using our general definition to quaternion gradients, substitute $h=0$ and $g = 3/4$ to compute the gradient of  
$f(\wz_{\text{RC}}(n))=e_{\text{RC}}(n)e_{\text{RC}}^*(n)$ and to obtain the update law
\begin{equation}
	\wz_{\text{RC}}(n+1)=\wz_{\text{RC}}(n)+\frac{3}{4}\mu_{\text{RC}} e_{\text{RC}}^{*}(n)\xz(n).
	\label{eq:update_wrc}
\end{equation}
Note that other values to $g$ could have been used, since the constant that appears in the second term of the right-hand side of \eqref{eq:update_wrc} 
can be absorbed by the step-size $\mu_{\text{RC}}$. We choose this value to make easier the comparison with other algorithms in the literature.
The algorithm is summarized in Table \ref{tab:rc_wl_iqlms}.
Since the additional term which appears in the RC-WL-QLMS algorithm due to the non-commutative quaternion multiplication 
vanishes in this new method, we expect it to be less costly to compute.
\begin{table}[htb]
 \centering
   \setlength{\arrayrulewidth}{1\arrayrulewidth}  
   \setlength{\belowcaptionskip}{10pt}  
   \caption{\it Proposed WL-QLMS algorithm with lower complexity}
   \begin{tabular}{|l|} 
      \hline
      $\hat{d}_{\text{RC}}(n)=\wz_{\text{RC}}^{H}(n)\xz(n)$ \\
      $e_{\text{RC}}(n)= d(n)-\hat{d}_{\text{RC}}(n)$\\
      $\wz_{\text{RC}}(n+1)=\wz_{\text{RC}}(n)+\frac{3}{4}\mu_{\text{RC}} e_{\text{RC}}^{*}(n)\xz(n)$ \\
      \hline
   \end{tabular}
   \label{tab:rc_wl_iqlms}
\end{table}

In the next section, we show that the proposed algorithm corresponds to the 4-Ch-LMS algorithm written in the quaternion domain, and that they 
have the same computational complexity.

\subsection[The RC-WL-iQLMS and the 4-Ch LMS algorithms]{Comparison with the 4-Ch-LMS algorithm}
\label{sec:comp_opt_LMS}

In order to prove the equivalence between the proposed algorithm and 4-Ch-LMS, 
we must rewrite the equations \eqref{eq:drc}, \eqref{eq:erc} and \eqref{eq:update_wrc} 
to separate the imaginary numbers $i$, $j$ and $k$ from the real terms. 

Define the extended vectors
\begin{align*}
  &\dz_{\text{ext}}(n) = \text{col}\left(d_\text{R}(n), \hspace{0.1 cm} d_\text{I}(n), \hspace{0.1 cm} d_\text{J}(n), \hspace{0.1 cm} d_\text{K}(n)\right),\\
  &\vz_{\text{ext}}'(n) = \text{col}\left(v_\text{R}(n), \hspace{0.1 cm} v_\text{I}(n), \hspace{0.1 cm} v_\text{J}(n), \hspace{0.1 cm} v_\text{K}(n)\right)
\end{align*}
and note that $d(n) = \tz^T\dz_\text{ext}(n)$ and $v(n) = \tz^T\vz_\text{ext}'(n)$, where
\begin{equation}
      \tz = \text{col}\left(\hspace{-0.1 cm}\begin{array}{cccc} 1, & i, & j, & k \end{array}\hspace{-0.1 cm}\right).
      \label{eq:t_tans}
\end{equation}
Assume that the $d(n)$ is modeled as given by
\begin{equation*}
    d(n) = \wz_{\text{o},\text{RC}}^H\xz(n) + v(n),
\end{equation*}
which can be expressed as  
\begin{equation*}
	\tz^T\dz_\text{ext}(n) =  \tz^T \left( \Wz_{{\text{o}},{\text{RC}}}^{T} \xz(n) +  \vz_\text{ext}'(n)  \right).
\end{equation*}
$\Wz_{{\text{o}},{\text{RC}}}$ is the $4N \times 4$ real matrix, given by
\begin{equation*}
  \Wz_{{\text{o}},{\text{RC}}} = \left[\wz_{{\text{o}},{\text{RC}}_\text{R}} \hspace{0.1 cm} -\wz_{{\text{o}},{\text{RC}}_\text{I}} \hspace{0.1 cm} -\wz_{{\text{o}},{\text{RC}}_\text{J}}\hspace{0.1 cm}  -\wz_{{\text{o}},{\text{RC}}_\text{K}}\right].
\end{equation*}
Define the extended vectors 
\begin{align}
  &\ez_{{\text{RC}}_{\text{ext}}}(n) = \text{col}\left(e_{{\text{RC}}_\text{R}}(n), \hspace{0.1 cm} e_{{\text{RC}}_\text{I}}(n),\hspace{0.1 cm}  e_{{\text{RC}}_\text{J}}(n), \hspace{0.1 cm} e_{{\text{RC}}_\text{K}}(n)\right),
   \label{eq:erro_ext} \\
  &\hat{\dz}_{{\text{RC}}_{\text{ext}}}(n) = \text{col}\hspace{-0.1 cm}\left(\hspace{-0.1 cm}\hat{d}_{{\text{RC}}_\text{R}}(n),\hspace{0.1 cm}  \hat{d}_{{\text{RC}}_\text{I}}(n),\hspace{0.1 cm}  \hat{d}_{{\text{RC}}_\text{J}}(n), \hspace{0.1 cm} \hat{d}_{{\text{RC}}_\text{K}}(n)\right).
  \label{eq:est_ext}
\end{align}
Using \eqref{eq:t_tans}, \eqref{eq:erro_ext} and \eqref{eq:est_ext} in eqs. \eqref{eq:drc}, \eqref{eq:erc} and \eqref{eq:update_wrc}, one gets
\begin{align}
   &\Wz(n+1) \tz^{*} =\left(\Wz(n) + \frac{3}{4}\mu_\text{RC} \xz(n)\ez_{{\text{RC}}_{\text{ext}}}^T(n)\right) \tz^* ,
  \label{eq:update_4} \\
	&\tz^T\ez_{{\text{RC}}_{\text{ext}}}(n) 
	      =  \tz^T \left(\dz_{\text{ext}}(n)- \hat{\dz}_{{\text{RC}}_{\text{ext}}}(n)\right) 
	\label{eq:error_4}
\end{align}
and
\begin{equation}
	\tz^T\hat{\dz}_{{\text{RC}}_{\text{ext}}}(n) = \tz^T\left(\Wz^{T}(n) \xz(n) \right),
	\label{eq:dest_4}
\end{equation}
where 
\begin{equation*}
  \Wz(n) = \left[\wz_{{\text{RC}}_\text{R}}(n) \hspace{0.1 cm}-\wz_{{\text{RC}}_\text{I}}(n) \hspace{0.1 cm}  -\wz_{{\text{RC}}_\text{J}}(n) \hspace{0.1 cm} -\wz_{{\text{RC}}_\text{K}}(n)\right].
\end{equation*}
Removing the multiplication by $\tz$ in \eqref{eq:update_4}, \eqref{eq:error_4} and \eqref{eq:dest_4}, we obtain the 4-Ch-LMS algorithm. Thus the proposed algorithm is a rewriting 
(in the quaternion domain) of the 4-Ch-QLMS algorithm, just as the RC-WL-LMS of \cite{fernandoISWCS2010} is a rewriting of the 2-Ch-LMS algorithm \cite{cap_vitor}.


\subsection[The RC-WL-iQLMS and the WL-iQLMS algorithms]{Comparing the new algorithm to WL-iQLMS}
\label{sec:comparison}

In this section, we show that the new approach and WL-iQLMS are related by a linear transformation, but the latter is more costly to compute.
 
The WL-iQLMS algorithm is presented in Table \ref{tab:wl_iqlms}.
\begin{table}[htb]
 \centering
   \setlength{\arrayrulewidth}{1\arrayrulewidth}  
   \setlength{\belowcaptionskip}{10pt}  
   \caption{\it WL-iQLMS algorithm}
   \begin{tabular}{|l|} 
      \hline
      $\hat{d}_{\text{iQ}}(n)=\wz_{\text{iQ}}^{H}(n)\qz_\text{WL}(n)$ \\
      $e_{\text{iQ}}(n)= d(n)-\hat{d}_{\text{iQ}}(n)$\\
      $\wz_{\text{iQ}}(n+1)=\wz_{\text{iQ}}(n)+\frac{3}{4}\mu_{\text{iQ}} e_{\text{iQ}}^{*}(n)\qz_\text{WL}(n)$ \\
      \hline
   \end{tabular}
   \label{tab:wl_iqlms}
\end{table}

To show the relation between the algorithms, define $\Fz$ as
\begin{equation*}
  \Fz = \left[\begin{array}{cccc}  
               1 & i & j & k\\
               1 & i & -j & -k\\
               1 & -i & j & -k\\
	       1 & -i & -j & k\\  
	      \end{array}\right] \otimes \Iz_{N},
\end{equation*}
and note that $\Fz \Fz^H/4 = \Fz^H \Fz/4 =\Iz_{4N}$.
One can show that 
\begin{equation}
  \qz_\text{WL}(n)=\Fz\xz(n).
  \label{eq:qa_qrc}
\end{equation}
The system of equations from which the optimum solution $\mathbf{w}_{o, \text{iQ}}$ is obtained is given
by
\begin{equation*}
  E\{\qz_\text{WL}(n)\qz_\text{WL}^H(n)\}\mathbf{w}_{\text{o},\text{iQ}} = E\{\qz_\text{WL}(n) d^*(n) \}.
\end{equation*}
Left-multiplying both sides by $\Fz^H/4$ and using \eqref{eq:qa_qrc}, one gets
\begin{equation*}
  \underbrace{\frac{\Fz^H}{4}\Fz}_{\Iz_{4N}} E\{\xz(n)\xz^T(n)\}\Fz^H \mathbf{w}_{\text{o},\text{iQ}} = \underbrace{\frac{\Fz^H}{4}\Fz}_{\Iz_{4N}} E\{\xz(n)d^*(n) \}.
\end{equation*}
Recognizing $E\{\xz(n)\xz^T(n)\}$ and $E\{\xz(n)d^*(n) \}$ as the correlation matrix and cross-correlation vector of the algorithm proposed in Table \ref{tab:rc_wl_iqlms}, then
\begin{equation}
    \wz_{\text{o}, \text{RC}} = \Fz^H \mathbf{w}_{\text{o},\text{iQ}} \;\; \text{or} \;\; \mathbf{w}_{\text{o},\text{iQ}} = \Fz\wz_{\text{o}, \text{RC}}/4.
    \label{eq:w_iqrc_opt}
\end{equation}
%
%
%
%
%
%
%
%
%
%
Based on \eqref{eq:w_iqrc_opt}, define
\begin{equation*}
  \wz_{\text{iQ}}(n) = \Fz \wz_{\text{RC}}(n)/4
\end{equation*}
to show that
\begin{equation*}
  \hat{d}_{\text{iQ}}(n)= \wz_{\text{iQ}}^{H}\frac{ \Fz^H\Fz}{4}(n) \qz_\text{WL}(n) =  \wz_{\text{RC}}^{H}(n) \xz(n)= \hat{d}_{\text{RC}}(n),
  \label{eq_dhat}
\end{equation*}
and notice that $e_{\text{iQ}}(n)=e_{\text{RC}}(n)$. Finally, left multiplying the update equation of WL-iQLMS by $\Fz^H$, we obtain
\begin{equation*}
  \underbrace{\Fz^H\wz_{\text{iQ}}(n+1)}_{\wz_{\text{RC}}(n+1)}=\underbrace{\Fz^H\wz_{\text{iQ}}(n)}_{\wz_{\text{RC}}(n)} + \frac{3}{4}\mu_{\text{iQ}}\Fz^H\underbrace{\qz_\text{WL}(n)}_{\Fz\xz(n)} \underbrace{e_{\text{iQ}}^{*}(n)}_{e_{\text{RC}}^{*}(n)}.
\end{equation*}
Defining $\mu_{\text{RC}}=4\mu_\text{iQ}$, we show that the technique proposed in Table \ref{tab:rc_wl_iqlms} and  WL-iQLMS are expected to 
have the same performance. 
Since the proposed algorithm uses $\mathbf{x}(n)$ as a real data vector, many quaternion-quaternion operations are replaced 
by real-quaternion operations (see Table \ref{tab:complexity}), reducing the computational cost. 
For this reason, the new approach is named as the RC-WL-iQLMS algorithm.

In Table \ref{tab:complexity}, we present the computational complexity of some quaternion LMS algorithms proposed in the literature. 
Notice that the WL algorithms are about 4 times more costly to implement than their SL counterparts, and that the complexity of the RC techniques
are similar to that of QLMS and iQLMS. Since RC-WL-iQLMS and iQLMS require the same number of computations, RC-WL-QLMS can be used 
as a low-cost alternative for both SL and WL scenarios. Finally, note that 
4-Ch-LMS and the RC-WL-iQLMS have the same computational complexity.

\section{Analysis of quaternion algorithms using real regressor vector}
\label{analysis}

In this section, we study the convergence of quaternion WL algorithms which use the real regressor vector $\mathbf{x}(n)$ 
(see eq. \eqref{eq:real_regressor}). Using our general description to the quaternion gradient, we obtain some results that can be applied to design 
the algorithms. We also obtain simple equations to compute the EMSE and the MSD.

\subsection{Designing $\mu$ to guarantee the convergence in the mean}
\label{sec:mean_analysis}



Similar to the analysis applied to study the LMS algorithm\cite{sayed2008adaptive}, one can use the maximum eigenvalue of
$\mathbf{G}_{\text{ext}}\mathbf{C}_{{\mathbf{x}}_{\text{ext}}}$ (see eq. \eqref{eq:mean_wl3}) to define bounds for the step-size which guarantee 
the convergence in the mean. From this approach, we obtain
\begin{equation*}
  0 < \mu < 2/(\text{max}\{(g+h), (g-h)\}\lambda_{\text{Max}}(\mathbf{C}_{\mathbf{x}})),
\end{equation*}
which is particularized to  
\begin{equation*}
  0 < \mu_{\text{RC}_{\text{Q}}}, \mu_\text{RC}  < 8/3\lambda_{\text{Max}}(\mathbf{C}_{\mathbf{x}}),
\end{equation*}
for RC-WL-QLMS and for RC-WL-iQLMS.

Both algorithms present the same bounds for the step-size, but they have different spreading factors, given by
\begin{align*}
    & SF_{\text{RC}_{\text{Q}}} = 3\lambda_{\text{Max}}(\mathbf{C}_{\mathbf{x}})/\lambda_{\text{Min}}(\mathbf{C}_{\mathbf{x}}) \text{ and}\\
    & SF_\text{RC} = \lambda_{\text{Max}}(\mathbf{C}_{\mathbf{x}})/\lambda_{\text{Min}}(\mathbf{C}_{\mathbf{x}}).
\end{align*}
Since the analysis was proposed for both correlated and uncorrelated input, one must expect RC-WL-iQLMS to converge faster than the RC-WL-QLMS. 
Our simulations in Section \ref{sec:simulations} confirm this.
\begin{table}[htb]
 \centering
   \setlength{\arrayrulewidth}{2\arrayrulewidth}  
   \caption[Computational complexity of the algorithms]{\it Computational complexity in terms of real operations per iteration 
      ($N$ is the length of the SL data vector)}
   \begin{tabular}{|c|c|c|} 
      \hline
      \textbf{Algorithm} & $+$ & $\times$ \\
      \hline
      QLMS & $ 48N$ & $ 48N + 9$\\
      \hline
      WL-QLMS &  $192N$ & $192N + 9$ \\
      \hline
      RC-WL-QLMS & $ 32N+4$ & $ 32N+8$ \\
      \hline
      iQLMS & $ 32N$ & $ 32N + 4$\\
      \hline
      WL-iQLMS & $ 128N$ & $ 128N + 4$\\
      \hline
      RC-WL-iQLMS & $ 32N$ & $ 32N + 4$\\
      \hline
      4-Ch-LMS & $ 32N$ & $ 32N + 4$\\      
      \hline
   \end{tabular}
   \label{tab:complexity}
\end{table}

\subsection{Mean-square analysis of real regressor quaternion algorithms}
\label{sec:mean_square_analysis}

In order to perform a general second-order analysis for real-regressor-vector quaternion algorithms, 
we drop all subscripts and use matrix $\Gz$ to account for any possible gradient results. 
We particularize our results for RC-WL-QLMS and RC-WL-iQLMS from the equations obtained.

To start the analysis, recall eq. \eqref{eq:qlms_gen1}. Assume that we are only considering WL-QLMS algorithms with real regressor vector to write
\begin{multline}
	\wz(n +1) = \wz(n) + \mu [(g+h)\mathbf{x}(n) e^{*}(n)\\ - 2h \imag{\{\mathbf{x}(n) e^{*}(n)}\}],
	\label{eq:w_real1}
\end{multline}
where $\wz(n)$ identifies the coefficients of a general algorithm with real-data vector. Define 
$\wz_{\text{o}}$ as the vector of optimum coefficients. Recall eq. \eqref{eq:w_tilde}.
Subtracting \eqref{eq:w_real1} from $\wz_\text{o}$, and using
\begin{equation*}
	e(n) = d(n) - (\wz_\text{o}^H \mathbf{x}(n) + v(n)),
\end{equation*}
we obtain
\begin{equation*}
	\tilde{\wz}(n +1) = \tilde{\wz}(n) - \mu [ \boldsymbol{\eta}(n) + \mathbf{u}(n)],
\end{equation*}
where we defined
\begin{align*}
   &\etaz(n)=\xz(n)\left((g + h) v^{*}(n) - 2h \imag \left\lbrace v^{*}(n)\right\rbrace\right) \text{ and} \\
    &\uz(n)=\xz(n)\xz^T(n)\left( (g + h)\tilde{\wz}(n) - 2h\imag\left\lbrace \tilde{\wz}(n) \right\rbrace\right).
\end{align*}
$\etaz(n)$ and $\uz(n)$ are quaternion vectors which depend on the noise and on $\tilde{\wz}(n)$, respectively.
Since eq. \eqref{eq:w_real1} contains quaternion entities, we use its extended version to perform the analysis.
We begin by calculating $\etaz_{\text{ext}}(n)$, which we divide in two parts, i.e.,
\begin{equation*}
      \etaz_{\text{ext}}(n) = \etaz_{1}(n) + \etaz_{2}(n).
\end{equation*}
$\etaz_{1}(n)$ is the extended version of $(g + h)\xz(n)v^{*}(n)$, that can be expressed as
\begin{equation}
    \etaz_{1}(n) = (g + h)\vz_{\text{ext}}(n) \otimes \xz(n),
     \label{eq:V1}
\end{equation}
where we define
\begin{equation*}
    \vz_{\text{ext}}(n) = \text{col}\left(v_\text{R}(n), \hspace{0.1 cm} -v_\text{I}(n), \hspace{0.1 cm} -v_\text{J}(n), \hspace{0.1 cm} -v_\text{K}(n)\right). 
\end{equation*}
$\etaz_{2}(n)$ is the extended version of $-2 h \xz(n)\imag\{v^{*}(n)\}$, i.e.,
\begin{equation}
    \etaz_{2}(n) = 2  h \hspace{0.05 cm}\text{col}\left(0, \hspace{0.1 cm} v_\text{I}(n), \hspace{0.1 cm} v_\text{J}(n), \hspace{0.1 cm} v_\text{K}(n)\right) \otimes \xz(n).
    \label{eq:V2}
\end{equation}
Define
\begin{equation}
    \Hz = \text{diag}((g+h), -(g-h), -(g-h), -(g-h))
    \label{eq:H2}
\end{equation}
and
\begin{equation}
  \Hz_{\text{ext}} = \Hz \otimes \Iz_{4N}.
  \label{eq:H}
\end{equation}
Finally, using \eqref{eq:V1}, \eqref{eq:V2} and \eqref{eq:H}, $\etaz_\text{ext}(n)$ is given by
\begin{equation*}
      \etaz_{\text{ext}}(n) = \Hz \vz_{\text{ext}}(n) \otimes \xz(n) = \Hz_{\text{ext}}( \vz_{\text{ext}}(n) \otimes \xz(n))
\end{equation*} 
Similarly to $\etaz_{\text{ext}}(n)$, $\uz_{\text{ext}}(n)$ is written as the sum of two terms, $\uz_{1}(n)$ and $\uz_{2}(n)$, which are the extended versions of 
$(g+h)\xz(n)\xz^T(n) \tilde{\wz}(n)$ and $ -2h\xz(n)\xz^T(n)\imag \{\tilde{\wz}(n)\}$. Using the extended version of $\tilde{\wz}(n)$, the first term is given by 
\begin{equation}
    \uz_{1}(n) = (g+h)\left( \Iz_{4} \otimes \xz(n)\xz^T(n) \right)\tilde{\wz}_{\text{ext}}(n),
    \label{eq:Q1}
\end{equation}
while $\uz_{2}(n)$ is written as
\begin{equation}
    \uz_{2}(n)= -2h \left(\text{diag}(0, \hspace{0.1 cm} 1, \hspace{0.1 cm} 1, \hspace{0.1 cm} 1) \otimes \xz(n)\xz^T(n) \right) \tilde{\wz}_{\text{ext}}(n)
    \label{eq:Q2}
\end{equation}
Adding up equations \eqref{eq:Q1} and \eqref{eq:Q2}, we obtain 
\begin{align*}
	\mathbf{u}_{\text{ext}}(n) =  \Gz_{\text{ext}} \left( \Iz_4 \otimes \xz(n)\xz^T(n) \right) \tilde{\wz}_\text{ext}(n).
\end{align*}
Considering the extended matrices, eq. \eqref{eq:w_real1} can be replaced in the analysis by
\begin{equation}
    \tilde{\wz}_{\text{ext}}(n+1)=  \tilde{\wz}_{\text{ext}}(n) - \mu \left[\etaz_{\text{ext}}(n) + \uz_{\text{ext}}(n) \right].
    \label{eq:rc_wl_qlms_3}
\end{equation}
Note that eq. \eqref{eq:rc_wl_qlms_3} deals only with real entities, which is fundamental in the following analysis.

Similarly to the traditional analysis of the LMS algorithm, multiply eq. \eqref{eq:rc_wl_qlms_3} by its transpose and take the expectation, to 
obtain
\begin{equation}
\begin{aligned}
    \underbrace{E\{\tilde{\wz}_{\text{ext}}(n+1)\tilde{\wz}^T_{\text{ext}}(n+1)\}}_{\Sz(n+1)} =  \underbrace{E\{\tilde{\wz}_{\text{ext}}(n)\tilde{\wz}^T_{\text{ext}}(n)\}}_{\Sz(n)}\\
	      - \mu \underbrace{E\{\tilde{\wz}_{\text{ext}}(n)\uz^T_{\text{ext}}(n)\}}_{\Am}   - \mu \underbrace{E\{\uz_{\text{ext}}(n)\tilde{\wz}^T_{\text{ext}}(n)\}}_{\Bm}\\
	      + \mu^2 \underbrace{E\{\uz_{\text{ext}}(n)\uz^T_{\text{ext}}(n)\}}_{\Cm}  + \mu^2 \underbrace{E\{\etaz_{\text{ext}}(n)\etaz^T_{\text{ext}}(n) \}}_{\Dm}\\
	       - \mu \underbrace{E\{\etaz_{\text{ext}}(n)\tilde{\wz}^T_{\text{ext}}(n)\}}_{\Em}  - \mu \underbrace{E\{\tilde{\wz}_{\text{ext}}(n)\etaz^T_{\text{ext}}(n)\}}_{\Fm}\\
	      + \mu^2 \underbrace{E\{\etaz_{\text{ext}}(n)\uz^T_{\text{ext}}(n)\}}_{\Gm} + \mu^2 \underbrace{E\{\uz_{\text{ext}}(n)\etaz^T_{\text{ext}}(n) \}}_{\Hm},
	    \label{eq:var}
\end{aligned}
\end{equation}
where we define $\Sz(n)=E\{\tilde{\wz}_{\text{ext}}(n)\tilde{\wz}^T_{\text{ext}}(n)\}$ to simplify the notation.
From eq. \eqref{eq:var} the second-order model for small step-sizes is derived. However, 
some approximations are required to proceed with the analysis, which are presented in the Assumption I next.

\textbf{Assumption I:} Assume that the sequence $\{\xz(n)\}$ is Gaussian, stationary and zero-mean, and that $\{\xz(n)\}$ and $\{v(n)\}$ 
are independent from each other. Additionally, assume that $E\{\vz_{\text{ext}}(n)\vz_{\text{ext}}^T(n)\}=\sigma_v^2 \Iz_4$ and that $\mu$ is small enough such that
$\mathbf{x}(n)$ and $\mathbf{w}(n)$ are approximately independent.

Based on Assumption I, the terms of eq. \eqref{eq:var} are studied as follows\footnote{Note that the we drop the time coefficients to simplify the notation.}:

\begin{enumerate}
 \item Term $\Am$ -- This term can be approximated as
       \begin{equation}
	  \begin{aligned}
	      \Am & = E\{ E\{ \tilde{\wz}_{\text{ext}}\tilde{\wz}^T_{\text{ext}}
			    \left( \Iz_4 \otimes \xz\xz^T \right)|\xz\} \}\Gz_{\text{ext}}\\
			& \approx  E\{ E\{\tilde{\wz}_{\text{ext}}\tilde{\wz}^T_{\text{ext}}\}
			    \left( \Iz_4 \otimes \xz\xz^T \right)\}\Gz_{\text{ext}}\\
			& \approx \Sz(n) E\{\left( \Iz_4 \otimes \xz(n)\xz^T(n) \right)\}\Gz_{\text{ext}}.
	  \end{aligned}
	  \label{eq:A}
       \end{equation}
       However, note that 
       \begin{equation*}
	  E\{ \Iz_4 \otimes \xz(n)\xz^T(n) \} = \Iz_4 \otimes E\{\xz(n)\xz^T(n)\}.
       \end{equation*}
       Using equations \eqref{eq:Cx} and \eqref{eq:Cx_ext}, eq. \eqref{eq:A} reduces to
       \begin{equation*}
	  \Am = \Sz(n) \Cz_{\xz_{\text{ext}}} \Gz_{\text{ext}}.
       \end{equation*}
  
  \item  Term $\Bm$ -- Similarly to $\Am$, term $\Bm$ reduces to 
	  \begin{equation*}
	     \Bm =  \Gz_\text{ext}\Cz_{{\xz}_{\text{ext}}}\Sz(n).
          \end{equation*}
          
  \item Term $\Cm$ -- This term can be rewritten as 
	\begin{equation}
	    \hspace{-0.4 cm}\begin{aligned}
	      \Cm & = \Gz_{\text{ext}} E\{ \left( \Iz_4 \otimes \xz\xz^T \right)
			    \tilde{\wz}_{\text{ext}}\tilde{\wz}^T_{\text{ext}}\left( \Iz_4 \otimes \xz\xz^T \right) \}\Gz_{\text{ext}}\\
			 & \cong \Gz_{\text{ext}} E\{ \left( \Iz_4 \otimes \xz\xz^T \right)
			\Sz(n)\left( \Iz_4 \otimes \xz\xz^T \right) \}\Gz_{\text{ext}}.
	    \end{aligned}
      \label{eq:C}
	\end{equation}
	
   \item Term $\Dm$ -- We can simplify $\Dm$ using
	 \begin{equation}
	    \begin{aligned}
		  \Dm & = E\{\etaz_{\text{ext}}(n)\etaz^T_{\text{ext}}(n) \}\\
		   & = \Hz_{\text{ext}} E\{\left( \vz_{\text{ext}} \otimes \xz\right)\left( \vz_{\text{ext}}^T \otimes \xz^T\right)\}\Hz_{\text{ext}}.
	    \end{aligned}
	    \label{eq:D}
	 \end{equation}
	 Using property \ref{property_3} of Kronecker products and the Assumption I, we simplify eq. \eqref{eq:D} to
	 \begin{equation*}
	    \begin{aligned}
		\Dm & = \Hz_{\text{ext}}\left( E\{\vz_{\text{ext}}\vz^T_{\text{ext}}\} \otimes E\{\xz\xz^T\}\right) \Hz_{\text{ext}}\\
			  & = \sigma_v^2\Hz_{\text{ext}} \left( \Iz_4 \otimes E\{\xz\xz^T\}\right) \Hz_{\text{ext}}\\
			  & = \sigma_v^2 \Hz_{\text{ext}}\Cz_{\xz_{\text{ext}}}\Hz_{\text{ext}}.
	     \end{aligned}
	 \end{equation*}
	 
     \item Term $\Em$ -- Using Kronecker product properties, we can rewrite $\Em$ as
	    \begin{equation*}
	      \begin{aligned}
		  \Em  = E\{\etaz_{\text{ext}}\tilde{\wz}^T_{\text{ext}}\}
		       = \Hz_{\text{ext}} E\{(\vz_{\text{ext}} \otimes \xz )\tilde{\wz}_{\text{ext}}^T \}.
	      \end{aligned}
	    \end{equation*}
	    From Assumption I, the elements of $\vz_{\text{ext}}$ are zero-mean random variables and are independent of $\xz$ and $\tilde{\wz}_{\text{ext}}$. 
	    In this case, $\Em$ is given by
	     \begin{equation*}
		  \begin{aligned}
		      \Em &  = \Hz_{\text{ext}}\left(  E\{\vz_{\text{ext}} \otimes \Iz_{N} \} E\{(\Iz_{4} \otimes\xz) \tilde{\wz}_{\text{ext}}^T \} \right)\\
			  &  = \Hz_{\text{ext}}(  E\{\vz_{\text{ext}} \} \otimes \Iz_{N}) E\{(\Iz_{4} \otimes\xz) \tilde{\wz}_{\text{ext}}^T \}
		  \end{aligned}
	    \end{equation*}
	    which results in a $16N \times 16N$ null matrix.
	    
      \item Terms $\Fm$, $\Gm$ and $\Hm$ -- The same argument is also applicable to these terms, which all result in $16N \times 16N$ null matrices.
	    
\end{enumerate}

After calculating all the terms of eq. \eqref{eq:var}, we can substitute the results to obtain
\begin{align}
  \Sz(n+1) \approx \Sz(n) - \mu \left[ \Sz(n)\Cz_{{\xz}_{\text{ext}}}\Gz_{\text{ext}} + \Gz_{\text{ext}}\Cz_{{\xz}_{\text{ext}}}\Sz(n)\right] \nonumber \\
			    +\mu^2\Cm + \mu^2\sigma_v^2 \Hz_{\text{ext}}\Cz_{{\xz}_{\text{ext}}}\Hz_{\text{ext}}.
  \label{eq:model}
\end{align}
Note that the three first terms on the right-hand side are linear in $\Sz(n)$. Assuming the small step-size condition \cite{sayed2008adaptive}, we can 
neglect the term $\mu^2\Cm$, which leads to the small step-size model
\begin{align}
   \Sz(n+1) \approx \Sz(n) - \mu \Sz(n)\Cz_{{\xz}_{\text{ext}}}\Gz_{\text{ext}} - \mu \Gz_{\text{ext}}\Cz_{{\xz}_{\text{ext}}}\Sz(n) \nonumber\\
			     + \mu^2\sigma_v^2 \Hz_{\text{ext}}\Cz_{{\xz}_{\text{ext}}}\Hz_{\text{ext}},
   \label{eq:model_small}
\end{align}
where the initialization corresponds to
$\Sz(0)=\tilde{\wz}_{\text{ext}}(0)\tilde{\wz}_{\text{ext}}^T(0)$.

Using \eqref{eq:model_small}, two theoretical quantities can be calculated at each time instant: the excess mean-square error (EMSE) \cite{sayed2008adaptive},
\begin{equation*}
  \zeta(n) = \text{Tr}(\Sz(n)\Cz_{{\xz}_{\text{ext}}}),
\end{equation*}
and the 
mean-square deviation (MSD) \cite{sayed2008adaptive}
\begin{equation*}
  \chi(n) = \text{Tr}(\Sz(n)).
\end{equation*}
From the small step-size model, we can also deduce the EMSE and MSD stead-state values. 
Assume that for $n \rightarrow \infty$, $\Sz(n+1) \approx \Sz(n) = \Sz(\infty)$, such that eq. \eqref{eq:model_small} reduces to
\begin{equation}
  \Sz(\infty)\Cz_{{\xz}_{\text{ext}}}\Gz_{\text{ext}} + \Gz_{\text{ext}}\Cz_{{\xz}_{\text{ext}}}\Sz(\infty)
			     = \mu\sigma_v^2 \Hz_{\text{ext}}\Cz_{{\xz}_{\text{ext}}}\Hz_{\text{ext}}.
  \label{eq:steady_state}
\end{equation}
Multiply eq. \eqref{eq:steady_state} by $\Gz_{\text{ext}}^{-1}$ from the right and take the trace. Using the trace property $\text{Tr}(\Az\Bz)=\text{Tr}(\Bz\Az)$, we obtain
\begin{equation}
  \zeta(\infty) \approx \mu\sigma_v^2 \text{Tr}(\Hz_{\text{ext}}\Cz_{{\xz}_{\text{ext}}}\Hz_{\text{ext}}\Gz_{\text{ext}}^{-1})/2.
  \label{eq:zeta}
\end{equation}
In eq. \eqref{eq:zeta}, use equations \eqref{eq:G_ext}, \eqref{eq:Cx_ext} and \eqref{eq:H2} to express
\begin{equation}
  \Hz_{\text{ext}}\Cz_{{\xz}_{\text{ext}}}\Hz_{\text{ext}}\Gz_{\text{ext}}^{-1} = (\Hz^2 \Gz^{-1}) \otimes \Cz_{\xz}.
  \label{eq:zeta1}
\end{equation}
Applying the trace property of the Kronecker product in eq. \eqref{eq:zeta1} and substituting the result in \eqref{eq:zeta}, one gets
\begin{equation*}
  \zeta(\infty) \approx \frac{\mu\sigma_v^2 \text{Tr}(\Hz^2\Gz^{-1})\text{Tr}(\Cz_{\xz})}{2} = (2g - h)\mu\sigma_v^2 \text{Tr}(\Cz_{\xz}).
\end{equation*}
Analogously, we compute the MSD by right multiplying both sides of eq.\eqref{eq:steady_state} by $\Gz^{-1}_{\text{ext}}\Cz_{{\xz}_{\text{ext}}}^{-1}$ and taking the trace, i.e.,
\begin{align*}
	  \chi(\infty) & \approx \mu \sigma_v^2\text{Tr}(\Hz_{\text{ext}}\Cz_{{\xz}_{\text{ext}}}\Hz_{\text{ext}}\Gz^{-1}_{\text{ext}}\Cz_{{\xz}_{\text{ext}}}^{-1})/2\\
		       & = \mu \sigma_v^2\text{Tr}(\Hz^2\Gz^{-1})\text{Tr}(\Cz_{{\xz}}\Cz_{{\xz}}^{-1})/2\\
		       & = 4N(2g-h)\mu \sigma_v^2 .
\end{align*}
Note that using this approach we obtain simple equations to calculate the EMSE and the MSD for small step-size, which depend only on $N$ and on the matrix $\Cz_{\xz}$, 
similarly to the LMS steady-state equations.

\subsection{Choosing the step-size}
\label{sec:step-size}

Recalling eq. \eqref{eq:model} and assuming that all variables are Gaussian, we can use properties of fourth-order Gaussian vectors 
\cite{sayed2008adaptive} to obtain an approximation for $\Cm$ (see eq.\eqref{eq:C}) and improve the accuracy of the proposed model. 
For this purpose, assume that $\xz(n)$ is a correlated Gaussian vector, and recall that the auto-correlation matrix $ \Cz_\mathbf{x} = E\{\mathbf{x}(n)\mathbf{x}^H(n)\}$ 
is symmetric and non-negative definite. This fact implies that there must exist an
unitary matrix $\Zz$, such that 
\begin{equation*}
  \Zz \Zz^T = \Zz^T \Zz = \Iz_{4N},
\end{equation*}
which diagonalizes $\Cz_{\xz}$, as given by
\begin{equation}
  \Cz_{\xz} = \Zz\Lbd \Zz^T ,
  \label{eq:R1}
\end{equation}
where $\Lbd$ is a diagonal matrix such that $\text{diag}(\Lbd)= [\lambda_1 \hspace{0.1 cm} \lambda_2 \hdots \lambda_{4N}]^T$, and $\lambda_i \geq 0$, $i=1, \hspace{0.1 cm} 2, \dots, 4N$
are the eigenvalues of $\Cz_{\xz}$. Defining $\xz'(n) = \Zz^T \xz(n)$, one can show 
\begin{equation*}
  E\{\xz'(n) \xz'^{T}(n)\}=E\{\Zz^T\xz(n)\xz^T(n)\Zz\} =\Zz^T \Cz_{\xz} \Zz = \Lbd.
\end{equation*} 
Since a linear transformation of a Gaussian vector is also Gaussian, $\xz'(n)$ is a Gaussian vector whose elements are independent from each other.
Define the extended matrix $\Zz_{\text{ext}} = \Iz_4 \otimes \Zz$. Multiplying $\Cm$ by $\Zz^T_{\text{ext}}$ on the left and by $\Zz_{\text{ext}}$ on the right, we obtain
\begin{align}
    &\Cm' = \Zz^T_{\text{ext}}\Cm\Zz_{\text{ext}} = \nonumber \\ 
    &=  \Zz^T_{\text{ext}}\Gz_{\text{ext}} E\{ \left( \Iz_4 \otimes \xz \xz^T \right)\Sz(n)\left( \Iz_4 \otimes \xz \xz^T \right) \}\Gz_{\text{ext}}\Zz_{\text{ext}}
    \label{eq:C_linha}
\end{align}
Defining $\Sz'(n)=\Zz^T\Sz(n)\Zz$ and noting that 
\begin{equation*}
  \Zz_{\text{ext}}^T\Gz_{\text{ext}}\Zz_{\text{ext}}=(\Iz_4 \otimes \Zz^T)(\Gz
	  \otimes \Iz_{4N})(\Iz_4 \otimes \Zz)=\Gz_{\text{ext}}
\end{equation*}
and
\begin{equation*}
      \begin{aligned}
	      \Zz^T_{\text{ext}}\left( \Iz_4 \otimes \xz\xz^T\right) \Zz_{\text{ext}} & = (\Iz_4 \otimes \Zz)^T(\Iz_4 \otimes \xz\xz^T)(\Iz_4 \otimes \Zz)\\
										  & = \Iz_4 \otimes \Zz^T\xz\xz^T\Zz\\
										  & =\Iz_4 \otimes \xz'\xz'^{T},
      \end{aligned}
\end{equation*}
we rewrite \eqref{eq:C_linha} as given by
\begin{equation}
  \Cm' = \Gz_{\text{ext}}E\left\lbrace  \left( \Iz_4 \otimes \xz'\xz'^{T} \right)\Sz'(n)\left( \Iz_4 \otimes \xz'\xz'^{T} \right) \right\rbrace \Gz_{\text{ext}}.
  \label{eq:C_dif}
\end{equation}
Each element of the matrix in the argument of \eqref{eq:C_dif} can be expressed as sums of terms
\begin{equation}
  E\{x'_{k1}(n) x'_{k2}(n) x'_{k3}(n) x'_{k4}(n)\}s'_{ij}(n),
  \label{eq:expc}
\end{equation}
where the $x'_{k}(n)$ represent the elements of $\xz'(n)$. $s'_{ij}(n)$ are the entries of $\Sz'(n)$. Since the elements of $\xz'(n)$ are 
independent and zero-mean, eq. \eqref{eq:expc} is different from zero if $k1=k2=k3=k4$, or $k1=k2$ and $k3=k4$, or $k1=k3$ and $k2=k4$, or $k1=k4$ and $k2=k3$. 
Using this result and eq. \eqref{eq:R1}, after some algebra manipulation, eq. \eqref{eq:C_dif} reduces to
\begin{equation*}
  \Cm' = \Gz_{\text{ext}} \left[ \Lbd_{\text{ext}} \text{Tr}\left( \Sz'(n)\Lbd_{\text{ext}}\right)
		      + 2\Lbd_{\text{ext}}\Sz'(n)\Lbd_{\text{ext}} \right]  \Gz_\text{ext}.
\end{equation*}
where $\Lbd_{\text{ext}} = \Iz_{4} \otimes \Lbd$. Recall eq. \eqref{eq:model} and multiply it by $\Zz^T_{\text{ext}}$ on the right and by $\Zz_{\text{ext}}$ on the left. Noting that
\begin{equation*}
  \Zz_{\text{ext}}^T\Hz_{\text{ext}}\Zz_{\text{ext}}=(\Iz_4 \otimes \Zz)^T(\Hz \otimes \Iz_{4N})(\Iz_4 \otimes \Zz)=\Hz_{\text{ext}},
\end{equation*}
we can write eq. \eqref{eq:model} as
\begin{multline}
      \Sz'(n+1) \approx  \Sz'(n) - \mu \left[ \Sz'(n)\Lbd_{\text{ext}}\Gz_{\text{ext}} + \Gz_{\text{ext}}\Lbd_{\text{ext}}\Sz'(n)\right]\\
		   	     + \mu^2\Gz_{\text{ext}}\left[\Lbd_{\text{ext}}\text{Tr}(\Sz'(n)\Lbd_{\text{ext}}) + 2\Lbd_{\text{ext}}\Sz'(n)\Lbd_{\text{ext}} \right]\Gz_{\text{ext}}\\
		   	     + \mu^2\sigma_v^2 \Hz_{\text{ext}}\Lbd_{\text{ext}}\Hz_{\text{ext}},
     \label{eq:model_mod}
\end{multline}
with initialization
\begin{equation*}
   \Sz'(0)=\Zz_{\text{ext}}^T\wz_{\text{ext}}(0)\wz_{\text{ext}}^T(0)\Zz_{\text{ext}}.
\end{equation*}
Taking only the diagonal $\sz'(n)=\text{diag}(\Sz'(n))$ in eq.\eqref{eq:model_mod}, we obtain a simplified recursion
\begin{equation}
\begin{aligned}
      \sz'(n+1) = \left[ \Iz_{16N} - 2\mu\Lbd_{\text{ext}}\Gz_{\text{ext}} +\mu^2\Gz^2_{\text{ext}}\boldsymbol{\ell}_{\text{ext}}\boldsymbol{\ell}_{\text{ext}}^T\right.  \\
			     \left. + 2 \mu^2\Gz^2_{\text{ext}}\Lbd_{\text{ext}}^2 \right] \sz'(n) + \mu^2\sigma_v^2\Hz^2_{\text{ext}}\boldsymbol{\ell}_{\text{ext}},
      \label{eq:model_mod1}
\end{aligned}
\end{equation}
where
\begin{equation*}
    \boldsymbol{\ell}_{\text{ext}} =\text{col}(1, \hspace{0.1 cm} 1, \hspace{0.1 cm} 1, \hspace{0.1 cm} 1) \otimes \text{diag}(\Lbd)
\end{equation*}
and $\sz'(0)=\text{diag}(\Sz'(0))$.

We can study the system matrix of eq. \eqref{eq:model_mod1}, i.e.,
\begin{align*}
  \boldsymbol{\Gamma}  = \Iz_{16N} - 2\mu\Lbd_{\text{ext}}\Gz_{\text{ext}} +\mu^2\Gz^2_{\text{ext}}\boldsymbol{\ell}_{\text{ext}}\boldsymbol{\ell}_{\text{ext}}^T,
			     + 2 \mu^2\Gz^2_{\text{ext}}\Lbd_{\text{ext}}^2\\
       = \left( \Iz_{16N} - \mu\Lbd_{\text{ext}}\Gz_{\text{ext}}\right)^2  +\mu^2\Gz^2_{\text{ext}}\boldsymbol{\ell}_{\text{ext}}\boldsymbol{\ell}_{\text{ext}}^T,
			     +  \mu^2\Gz^2_{\text{ext}}\Lbd_{\text{ext}}^2
\end{align*}
to define a bound for the step-size which guarantees the stability in the variance. Using the $\ell_1$-norm \cite{meyerMAALA00}, we can find an upper 
bound for the largest eigenvalue  $\nu_l$ of $\boldsymbol{\Gamma}$, i.e.,
\begin{equation*}
    \max_{1 \leq l \leq 16N} |\nu_l| \leq ||\boldsymbol{\Gamma}||_1= \max_{1 \leq l \leq 16N} \sum_{\text{m}=1}^{16N}|\gamma_{l\text{m}}|,
\end{equation*}
where $\gamma_{l\text{m}}$ are the elements of $\boldsymbol{\Gamma}$. A conservative range of values for $\mu$, which guarantee the stability, requires 
that $||\boldsymbol{\Gamma}||_1 \leq 1$. 
Observe that the $l$-th column of $\boldsymbol{\Gamma}$ has entries $\mu^2 g_{\text{ext}_{l,l}}^2\lambda_{\text{ext}_{l}}\lambda_{\text{ext}_{m}}$, if $l \neq m$, and 
$(1-\mu g_{\text{ext}_{l,l}}\lambda_{\text{ext}_{l}})^2+ \mu^2 g_{\text{ext}_{l,l}}^2\lambda_{\text{ext}_{l}}^2 + \mu^2g_{\text{ext}_{l,l}}^2\lambda_{\text{ext}_{l}} \sum_{m=1}^{16N}\lambda_{\text{ext}_{m}}$, if $l = m$, where 
$g_{\text{ext}_{l,l}}$ are the diagonal elements of $\Gz_{\text{ext}}$. Note that
\begin{multline*}
     \hspace{-0.2 cm}\sum_{\text{m}=1}^{16N}\hspace{-0.12 cm}|\gamma_{l\text{m}}|\hspace{-0.1 cm}
     =\hspace{-0.1 cm}(1\hspace{-0.05 cm}-\hspace{-0.05 cm}\mu g_{\text{ext}_{l,l}}\lambda_{\text{ext}_{l}})^2 \hspace{-0.07 cm}+ \hspace{-0.07 cm}\mu^2  g_{\text{ext}_{l,l}}^2\hspace{-0.15 cm}\left(\hspace{-0.05 cm}\lambda_{\text{ext}_{l}}^2\hspace{-0.1 cm} + \hspace{-0.05 cm}\lambda_{\text{ext}_{l}} \hspace{-0.15 cm}\sum_{\text{m}=1}^{16N}\hspace{-0.1 cm}\lambda_{\text{ext}_{\text{m}}}\hspace{-0.1 cm}\right)\\
     =(1-\mu g_{\text{ext}_{l,l}}\lambda_{\text{ext}_{l}})^2+ \mu^2 g_{\text{ext}_{l,l}}^2\left( \lambda_{\text{ext}_{l}}^2 + \lambda_{\text{ext}_{l}} \text{Tr}(\boldsymbol{\Lambda}_{\text{ext}})\right).
\end{multline*}
Thus, the recursion in eq. \eqref{eq:model_mod1} is stable if
\begin{equation*}
    (1-\mu g_{\text{ext}_{l,l}}\lambda_{\text{ext}_{l}})^2+ \mu^2 g_{\text{ext}_{l,l}}^2\left( \lambda_{\text{ext}_{l}}^2 + \lambda_{\text{ext}_{l}} \text{Tr}(\Lbd_{\text{ext}})\right) \leq 1,
\end{equation*}
for $1 \leq l \leq 16N$. 
Recalling that 
\begin{equation*}
    \text{Tr}(\Cz_{{\xz}_{\text{ext}}})=\text{Tr}(\Lbd_{\text{ext}})=\text{Tr}(\Iz_4 \otimes \Cz_{\xz} )=4\text{Tr}(\Cz_{\xz})
\end{equation*}
and after some manipulation, the condition simplifies to
\begin{equation*}
    \begin{array}{cc}
	\mu \leq 1/g_{\text{ext}_{l,l}}(\lambda_{\text{ext}_{l}} + 2\text{Tr}(\Cz_{\xz})), & 1 \leq l \leq 16N.
    \end{array}
\end{equation*}
The smallest bound occurs for $g_{\text{ext}_{l,l}}=\text{max}\{(g+h), (g-h)\}$ and $\lambda_{\text{ext}_{l}}=\lambda_{{\text{ext}}_\text{max}}$:
\begin{equation}
	\mu \leq 1/(\text{max}\{(g+h), (g-h)\}(\lambda_{{\text{ext}}_{\text{max}}} + 2\text{Tr}(\Cz_{\xz}))).
    \label{eq:step_size}
\end{equation}
Replacing $\lambda_{{\text{ext}}_{\text{max}}}$ by $\text{Tr}(\Cz_{\xz})$ we obtain a simpler but more conservative condition for stability,
\begin{equation*}
	0 < \mu < 1/(3\text{max}\{(g+h), (g-h)\}\text{Tr}(\Cz_{\xz})),
\end{equation*}
which guarantees stability in the variance. For RC-WL-QLMS and RC-WL-iQLMS, the step-size selection must respect
\begin{equation}
    \begin{array}{cc}
	0 < \mu_\text{iQ}, \mu_\text{RC} <4/(9\text{Tr}(\Cz_{\xz})),
    \end{array}
    \label{eq:murc1}
\end{equation}
since $\text{max}\{(g+h), (g-h)\}=3/4$ for both algorithms.

Note that analysis proposed in this section is valid for uncorrelated and correlated input. 
The results presented here can be extended to other real-regressor quaternion algorithms.

\section{Simulations}
\label{sec:simulations}

In order to compare the algorithms and show the accuracy of the proposed model,
we performed some simulations using $\nq$-improper processes. For this purpose, we conveniently define the elements of $\qz(n)$ as given by
\begin{align*}
    \qz_\text{R}(n) &= \boldsymbol{\rho}_1(n) + 0.1\boldsymbol{\rho}_2(n) + 0.2\boldsymbol{\rho}_3(n) + 0.3\boldsymbol{\rho}_4(n)\\
    \qz_\text{I}(n) &= 0.1\boldsymbol{\rho}_1(n) + \boldsymbol{\rho}_2(n) + 0.1\boldsymbol{\rho}_3(n) + 0.1\boldsymbol{\rho}_4(n)\\
    \qz_\text{J}(n) &= 0.2\boldsymbol{\rho}_1(n) + 0.1\boldsymbol{\rho}_2(n) + \boldsymbol{\rho}_3(n) + 0.1\boldsymbol{\rho}_4(n)\\
    \qz_\text{K}(n) &= 0.3\boldsymbol{\rho}_1(n) + 0.1\boldsymbol{\rho}_2(n) + 0.1\boldsymbol{\rho}_3(n) + \boldsymbol{\rho}_4(n), 
\end{align*}
where $\boldsymbol{\rho}_l(n)$ are $4 \times 1$ vectors. We assume that the elements of each $\boldsymbol{\rho}_l(n)$ are zero-mean, Gaussian, 
i.i.d. and with unitary variance, and that $\boldsymbol{\rho}_l(n)$ and $\boldsymbol{\rho}_m(n)$ are independent from each other 
$\forall l \neq m$. That means that the components of each quaternion in $\qz(n)$ are correlated, but different quaternions are uncorrelated.
The desired sequence $d(n)$ and $\qz(n)$ are jointly-$\nq$-improper processes, since $d(n)$ is obtained with
\begin{equation*}
    d(n) = \wz_o^H\qz_\text{WL}(n) + v(n).
\end{equation*}
$\wz_o$ is a $16 \times 1$ quaternion vector, where the elements of the first four quaternions are obtained from a $\text{Uniform}(0,1)$ distribution, and 
the elements of the other quaternions are obtained from $\text{Uniform}(0,10^{-2})$ distribution. Note that with this approach, we guarantee a WL
$d(n)$ sequence, but we give more emphasis to the SL contribution.  

The quaternion elements of $v(n)$ are zero-mean, Gaussian and i.i.d, with equal variance $\sigma_v^2=0.001/4$. We perform $200$ simulations to observe the EMSE and the MSD.
We compare QLMS, WL-QLMS, WL-iQLMS, RC-WL-QLMS, RC-WL-iQLMS and 4-Ch-LMS, and we plot our model for RC-WL-iQLMS and RC-WL-QLMS.
We adjust the WL-iQLMS and the QLMS algorithms to have the same convergence rate, and adjust the other WL-algorithms (and the 4-Ch-LMS)
to achieve the same steady-state EMSE. Table \ref{tab:mu} shows the step-sizes used in our simulations.
\begin{table*}[htb]
 \centering
   \setlength{\arrayrulewidth}{2\arrayrulewidth}  
   \setlength{\belowcaptionskip}{10pt}  
   \caption{\it Step-sizes used in the simulations ($\mu=10^{-5}$)}
   \begin{tabular}{|c|c|c|c|c|c|c|} 
      \hline
      \textbf{Algorithm} & QLMS & WL-QLMS & WL-iQLMS & RC-WL-QLMS & RC-WL-iQLMS & 4-Ch-LMS\\
      \hline
      \textbf{Step-size} & $\mu$ & $0.8\mu$ & $2\mu$ & $3.2\mu$ & $8\mu$ & $6\mu$\\    
      \hline
   \end{tabular}
   \label{tab:mu}
\end{table*}
\begin{figure}[hptb]
  \psfrag{EMSE (dB)}{EMSE(dB)}
  \psfrag{Iterations}{Iterations}
  \psfrag{QLMS}{\begin{scriptsize} QLMS \end{scriptsize}}
  \psfrag{WL-QLMS}{\begin{scriptsize} WL-QLMS \end{scriptsize}}
  \psfrag{RC-WL-QLMS}{\begin{scriptsize} RC-WL-QLMS \end{scriptsize}}
  \psfrag{RC-WL-iQLMS}{\begin{scriptsize} RC-WL-iQLMS \end{scriptsize}}
  \psfrag{WL-iQLMS}{\begin{scriptsize} WL-iQLMS \end{scriptsize}}
  \psfrag{4-Ch-LMS}{\begin{scriptsize} 4-Ch-LMS \end{scriptsize}}
  \psfrag{Modelo do RC-WL-iQLMS}{\begin{scriptsize} RC-WL-iQLMS Model \end{scriptsize}}
  \psfrag{Modelo do RC-WL-QLMS}{\begin{scriptsize} RC-WL-QLMS Model \end{scriptsize}}
  \includegraphics[width=\linewidth]{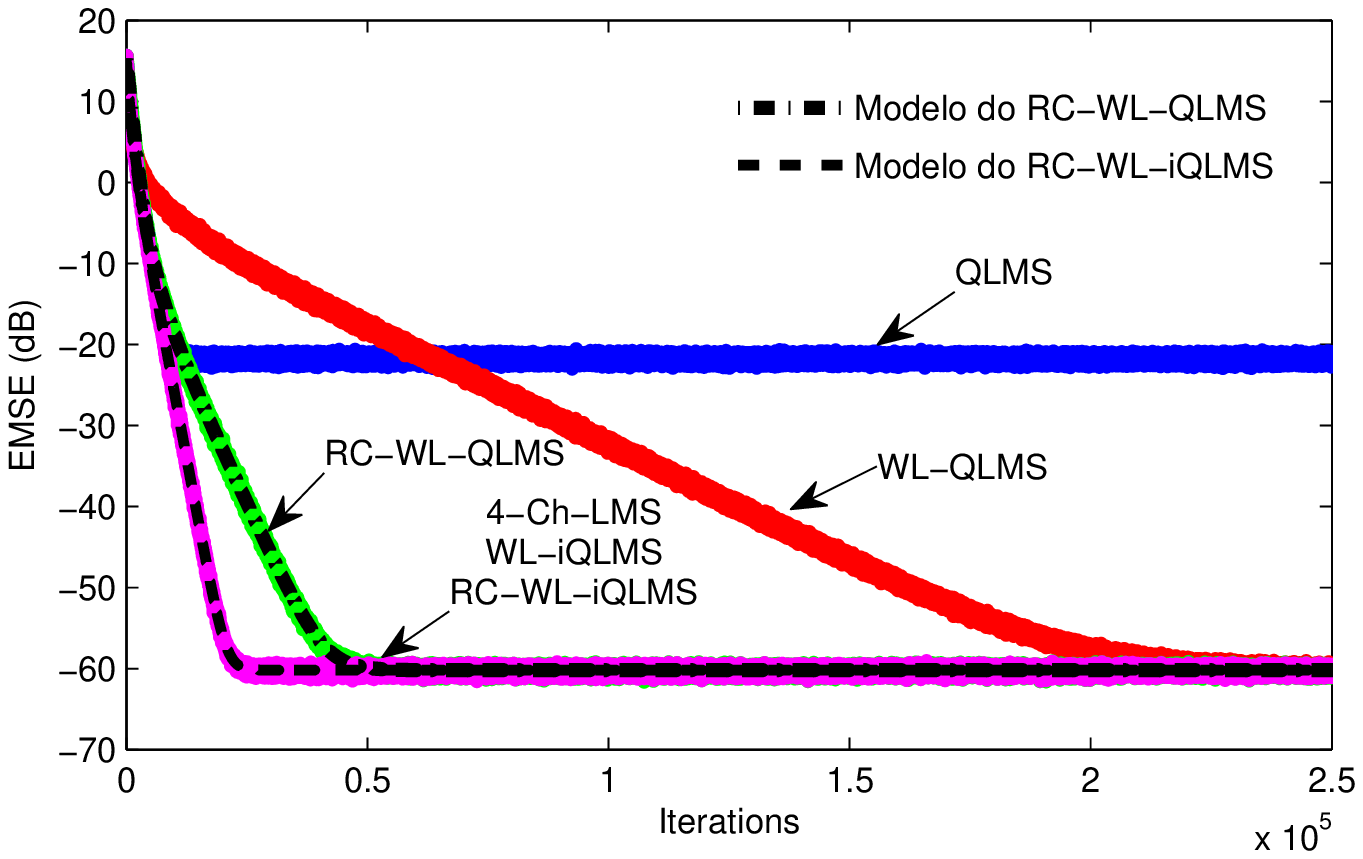}
  \caption{EMSE performance comparing QLMS, WL-QLMS, RC-WL-QLMS, RC-WL-iQLMS, WL -iQLMS , 4-Ch-LMS and the proposed second-order model  --Average of $200$ simulations.}
  \label{fig:EMSE}
\end{figure}
 \begin{figure}[hptb]
   \psfrag{MSD (dB)}{MSD(dB)}
   \psfrag{Iterations}{Iterations}
   \psfrag{QLMS}{\begin{scriptsize} QLMS \end{scriptsize}}
   \psfrag{WL-QLMS}{\begin{scriptsize} WL-QLMS \end{scriptsize}}
   \psfrag{RC-WL-QLMS}{\begin{scriptsize} RC-WL-QLMS \end{scriptsize}}
   \psfrag{RC-WL-iQLMS}{\begin{scriptsize} RC-WL-iQLMS \end{scriptsize}}
   \psfrag{WL-iQLMS}{\begin{scriptsize} WL-iQLMS \end{scriptsize}}
   \psfrag{4-Ch-LMS}{\begin{scriptsize} 4-Ch-LMS \end{scriptsize}}
   \psfrag{Modelo do RC-WL-iQLMS}{\begin{scriptsize} RC-WL-iQLMS Model \end{scriptsize}}
   \psfrag{Modelo do RC-WL-QLMS}{\begin{scriptsize} RC-WL-QLMS Model \end{scriptsize}}
   \includegraphics[width=\linewidth]{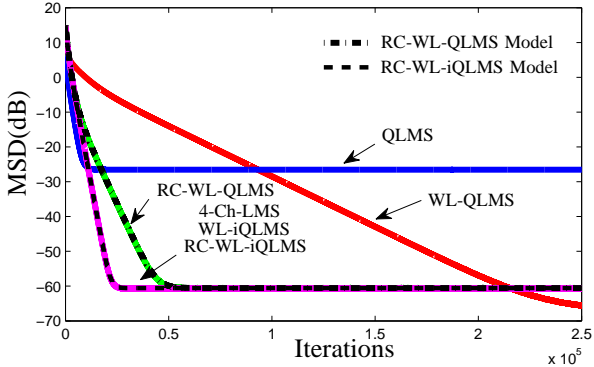}
   \caption{MSD performance comparing QLMS, WL-QLMS, RC-WL-QLMS, RC-WL-iQLMS, WL-iQLMS , 4-Ch-LMS and the proposed second-order model --Average of $200$ simulations.}
   \label{fig:MSD}
 \end{figure}

Note that the WL algorithms achieve lower EMSE and MSD performances than QLMS, and that WL-iQLMS, 4-Ch-LMS and RC-WL-iQLMS have the 
same converge rate, which is faster than WL-QLMS and RC-WL-QLMS. In addition, the proposed model (presented in the figures as a 
black dashed line) is accurate to describe the performance behavior of both the RC-WL-QLMS and the RC-WL-iQLMS algorithms.

\section{Conclusions}
\label{sec:conclusions}

In this paper, we developed a general representation to the quaternion gradients proposed in the literature. 
Using this approach, we proved that different gradients can be used to obtain the same algorithm.

We showed that the class of gradients from which the i-gradient takes part provides the fastest-converging WL algorithms when the 
correlation matrix has entries with only the real and one imaginary part, and when the WL regressor vector is real 
and obtained by the concatenation of the real and imaginary elements of the original SL data vector. 
The general gradient was applied to  devised the fastest-converging WL-QLMS algorithm with 
real-regressor vector -- the RC-WL-iQLMS algorithm -- and we showed that the proposed method corresponds to the 4-Ch-LMS written in the quaternion 
field. It was also shown that new algorithm corresponds to a lower-complexity version of WL-iQLMS, 4 times less 
costly to implement, and with the same complexity of the 4-Ch-LMS algorithm. 
In Section \ref{sec:mean_square_analysis} we extended the second-order analysis of \cite{fernandoISWCS2012} to any WL-QLMS algorithm using real 
data vector. This approach led to simple equations to compute largely applied figures of merit used in adaptive filtering, such as the EMSE and 
the MSD. The comparison between the model and the performance of the algorithms proved that the model is accurate.

\appendices

\section{Conditions to guarantee the positive semi-definiteness of $\mathbf{G}_{\mathrm{ext}}\mathbf{C}_{\mathrm{ext}}$}
\label{ap:GC}

When we use the real extended entities to obtain eq. \eqref{eq:mean_wl2} from \eqref{eq:equation_quat},
our goal is the application of mathematical tools from the real field to access the eigenvalue spread of $\mathbf{G}_{\text{ext}}\mathbf{C}_{\text{ext}}$, such that we can 
study the convergence of quaternion algorithms. For this purpose, in this appendix we show that $\mathbf{C}_{\text{ext}}$ is a 
positive semi-definite matrix. Then, we use this to prove that the diagonal entries of $\mathbf{G}_{\text{ext}}$ must be non-negative to 
make $\mathbf{G}_{\text{ext}}\mathbf{C}_{\text{ext}}$ positive semi-definite always. 

Initially, recall that the original quaternion correlation matrix is given by $E\{\mathbf{q}_{\text{WL}}(n)\mathbf{q}_{\text{WL}}^H(n)\}$. It is easy to notice that for any quaternion vector 
$\boldsymbol{\chi}$ with the proper dimension,
\begin{equation}
    \boldsymbol{\chi}^H E\{\mathbf{q}_{\text{WL}}(n)\mathbf{q}_{\text{WL}}^H(n)\} \boldsymbol{\chi} = 
    E\{|\mathbf{q}_{\text{WL}}^H(n)\boldsymbol{\chi}|^2\} \geq 0, \;\; \forall \boldsymbol{\chi},
    \label{eq:pos_def}
\end{equation}
such that the correlation matrix is positive semi-definite. 

From eq. \eqref{eq:pos_def}, one would expect $\mathbf{C}_{\text{ext}}$ also positive semi-definite, since it is the extended version 
of the correlation matrix. However, when we check the block structure of $\mathbf{C}_{\text{ext}}$, it is not trivial to show that
this property holds, since the extended matrix is not easily expressed as the product of a real vector by its transpose. 
To show that $\Cz_\text{ext}$ is also positive semi-definite, one can explicit the terms which appear in the computation of \eqref{eq:pos_def}, 
and then reorganize them to reveal
\begin{equation*}
   \boldsymbol{\chi}_\text{ext}^T \mathbf{C}_\text{ext}\boldsymbol{\chi}_\text{ext} =\boldsymbol{\chi}^H E\{\mathbf{q}_{\text{WL}}(n)\mathbf{q}_{\text{WL}}^H(n)\} \boldsymbol{\chi} \geq 0,
\end{equation*}
 where
$\boldsymbol{\chi}_\text{ext} = \text{col}(\hspace{-0.1 cm}\begin{array}{cccc} \boldsymbol{\chi}_\text{R}, & \boldsymbol{\chi}_\text{I}, & \boldsymbol{\chi}_\text{J}, & \boldsymbol{\chi}_\text{K} \end{array} \hspace{-0.1 cm})$.
Since  $\forall \boldsymbol{\chi}_\text{ext}$ can be obtained from the mapping of $\boldsymbol{\chi}$ to the real field, $\Cz_\text{ext}$ is also
positive semi-definite.

To define the cases on which $\mathbf{G}_{\text{ext}}\mathbf{C}_{\text{ext}}$ leads to stable \eqref{eq:mean_wl2}, we must find values for which 
the diagonal entries of $\mathbf{G}_{\text{ext}}$ make the product of matrices positive semi-definite. First, recall that using an unitary transformation 
\cite{meyerMAALA00}, one can show that $\mathbf{G}_{\text{ext}}\mathbf{C}_{\text{ext}}$ has the same eigenvalues of 
\begin{equation}
  \mathbf{G}_{\text{ext}}^{1/2}\mathbf{C}_{\text{ext}}\mathbf{G}_{\text{ext}}^{1/2}. 
  \label{eq:G_C_G}
\end{equation}
In this case, we just need to define the conditions on which \eqref{eq:G_C_G} is positive semi-definite to show that the elements of $\mathbf{G}_{\text{ext}}$ should be non-negative.
For this purpose, we use vector $\boldsymbol{\chi}_\text{ext}$ to check on which cases
\begin{equation*}
	\boldsymbol{\chi}_{\text{ext}}^H\mathbf{G}_{\text{ext}}^{1/2}\mathbf{C}_{\text{ext}}\mathbf{G}_{\text{ext}}^{1/2}\boldsymbol{\chi}_{\text{ext}} \geq 0.
\end{equation*}
We start showing that if $\mathbf{G}_{\text{ext}}$ has at least one negative diagonal entry, it is possible to find one vector for which 
$\mathbf{G}_{\text{ext}}^{1/2}\mathbf{C}_{\text{ext}}\mathbf{G}_{\text{ext}}^{1/2}$ is not positive semi-definite. Assume that only the $k$-th entry 
$\mathbf{G}_{\text{ext}}$ is $g_{\text{ext}_{k,k}} < 0$, such that
\begin{equation*}
    \mathbf{G}_{\text{ext}}^{1/2} \hspace{-0.1 cm}= \text{diag}(\hspace{-0.2 cm}\begin{array}{ccccc} \sqrt{g_{\text{ext}_{1,1}}},
    & \hspace{-0.2 cm}\hdots, & \hspace{-0.2 cm}j\sqrt{|g_{\text{ext}_{k,k}}|}, & \hspace{-0.2 cm}\hdots, & \hspace{-0.2 cm}\sqrt{g_{\text{ext}_{16N,16N}}} \end{array}),
\end{equation*}
where $j = \sqrt{-1}$. Using $\boldsymbol{\chi}_\text{ext} = \boldsymbol{\epsilon}_k$, where $\epsilon_k=1$ and the other entries are set to zero, one can easily show
\begin{equation*}
  \boldsymbol{\epsilon}_k^T\mathbf{G}_{\text{ext}}^{1/2}\mathbf{C}_{\text{ext}}\mathbf{G}_{\text{ext}}^{1/2}\boldsymbol{\epsilon}_k = g_{\text{ext}_{k,k}} c_{\text{ext}_{k, k}} < 0 
\end{equation*}
so that \eqref{eq:G_C_G} is not positive semi-definite. When $\mathbf{G}_{\text{ext}}$ has more negative diagonal entries, the proof that 
\eqref{eq:G_C_G} is non-positive is straightforward.
In this case, all $g_{\text{ext}_{k,k}}$ must be non negative, $1 \leq k \leq 16N$, to guarantee the stability of eq. \eqref{eq:mean_wl2}. The condition for this
is $(g + h) \geq 0$ and $(g-h) \geq 0$.      

%
%
%
%

\bibliographystyle{IEEEtran}
\bibliography{IEEEabrv,refs}

\end{document}